\documentclass[10pt,authoryear]{article}
\usepackage{appendix}
\usepackage{setspace}
\usepackage[font=small,  margin=2cm]{caption}
\usepackage[tbtags]{amsmath}
\usepackage{amsthm,amssymb,amsfonts}
\usepackage{natbib}
\usepackage{enumitem}
\usepackage{multirow}
\usepackage{lscape}
\usepackage{longtable}
\usepackage{subcaption}
\usepackage{caption}
\usepackage{graphicx}
\usepackage{epstopdf}
\usepackage{epsfig}
\usepackage{makecell}
\usepackage{exscale}
\usepackage{booktabs}
\usepackage{array}
\usepackage{fullpage}
\usepackage{url}
\usepackage{algorithm}
\usepackage{algorithmic}

\usepackage{bm}
\usepackage{smile}
\usepackage{mathtools}
\usepackage{wrapfig}
\usepackage{lipsum}
\usepackage{mathrsfs}
\usepackage{relsize}
\usepackage{dsfont}
\usepackage{multirow}
\usepackage{apalike}
\usepackage{chngcntr}
\usepackage{graphicx}
\usepackage{hyperref}
\usepackage[usenames,dvipsnames,svgnames,table]{xcolor}

\theoremstyle{remark}
\newtheorem{remark}{Remark}
\newtheorem{assumption}{Assumption}

\numberwithin{equation}{section}
\numberwithin{theorem}{section}
\numberwithin{corollary}{section}
\counterwithout{asmp}{section}
\numberwithin{definition}{section}

\usepackage{amsmath}
\allowdisplaybreaks[4]

\begin{document}

	\title{TransPCA for Large-dimensional  Factor Analysis with Weak Factors: Power Enhancement via Knowledge Transfer}

	\author{Yong He\footnotemark[1],~~Dong Liu\footnotemark[2],~~Yunjing Sun\footnotemark[3],~~Yalin Wang\footnotemark[3]}
\renewcommand{\thefootnote}{\fnsymbol{footnote}}
\footnotetext[1]{Institute for Financial Studies, Shandong University, China}
\footnotetext[2]{Nanyang Technological University, Singapore}	
\footnotetext[3]{Corresponding author. School of Mathematics, Shandong University, China; e-mail: {\tt wangyalin@mail.sdu.edu.cn}}

	\date{}
	\maketitle
Early work established convergence of the principal component
estimators of the factors and loadings up to a rotation for large dimensional
approximate factor models with weak factors in that the factor loading $\bLambda^{(0)}$ scales
sublinearly in the number $N$ of cross-section units, i.e., $\bLambda^{(0)\top}\bLambda^{(0)}/N^{\alpha}$  is positive definite
in the limit for some $\alpha\in (0,1)$. However, the established convergence rates for weak factors can be much slower especially for small $\alpha$.
This article proposes a Transfer Principal Component Analysis (TransPCA) method  for enhancing the convergence rates for weak factors by transferring knowledge from large number of available informative panel datasets, which should not be turned a blind eye on in this big data era. 
We aggregate useful information by analyzing a weighted average projection matrix of the estimated loading spaces from all informative datasets which is highly flexible and computationally efficient. Theoretically, we derive the convergence rates of the estimators of weak/strong loading spaces and factor scores. The results indicate that as long as  the auxiliary datasets are similar enough to the target dataset and the auxiliary sample size is sufficiently large, TransPCA estimators can achieve faster convergence rates in contrast to performing PCA solely on the target dataset. We also propose a Transfer Eigenvalue Difference (TransED) method to determine the number of weak factors based
on the maximum eigengap and two kind of methods for estimating the factor strengths, and  the consistency of the resultant estimators is also derived. To avoid negative transfer, we also investigate the
case that the informative datasets are unknown and provide a criterion  for selecting useful
datasets. Thorough simulation studies and {empirical analysis on real datasets in areas of  macroeconomic and finance}  are conducted to illustrate the usefulness of our proposed methods where large number of source panel datasets are naturally available.

\vspace{0.5em}
\section{Introduction}

    Large-dimensional approximate factor model has been well studied in the last two decades, which acts as a powerful tool for dimension reduction and plays an extremely important role in the era of big data nowadays. Since the seminal works of \cite{bai2002determining} and \cite{stock2002forecasting}, it has been witnessed that a large body of research focuses on
estimating the latent common variations in large panels in which the $N$ units observed over $T$ periods are cross-sectionally correlated, see for example,  \cite{bai2003inferential}, \cite{onatski2009testing}, \cite{ahn2013eigenvalue}, \cite{fan2013large}, \cite{bai2012statistical}, \cite{trapani2018randomized},  \cite{barigozzi2020sequential} and \cite{kong2024staleness}. In view of the typical heavy-tailedness characteristic of financial dataset,  there have been flourishing works in the past few years  on relaxing the moment conditions for factor analysis, see for example the endeavors by \cite{yu2019robust}, \cite{chen2021quantile}, \cite{he2022large}, \cite{barigozzi2024tail} and \cite{he2025huber}.

\subsection{ Closely Related Work}
Among the extensive studies on the estimation of factor models, few  pay attention to the availability of large number of auxiliary source datasets which could be informative for the target factor analysis.
 Actually, it is growing important to extract information and transfer knowledge from diversified data sources to embrace the big data era. Turning a blind eye to these informative datasets would lead to inefficiency in downstream statistical tasks such as prediction. To the best of our knowledge,  only a few works delve into  factor modeling with additional auxiliary sources datasets. For instance, the constrained factor model in \cite{tsai2010constrained} exerts complete or partial linear constraints on the factor loadings, according to the additional industrial sectors information of 10 U.S. companies in their real data analysis; the semi-parametric factor structure in \cite{connor2012efficient} and \cite{fan2016projected}  models the factor loadings as unknown functions of certain observable covariates, which utilizes auxiliary datasets  such as the market capitalization and price-earning ratio; \cite{yu2020network}  proposed a  network-linked  factor model   to incorporate auxiliary network information among the large-scale variables; \cite{duan2024target} developed a target PCA method to do factor analysis for a large target panel
with missing observations by extracting information from auxiliary panel datasets.
 Indeed, taking advantage of auxiliary information/datasets in factor modeling is still in its infancy and certainly deserves more attention. Meanwhile, the improper use of uninformative datasets should be avoided, which makes it even more challenging to extract knowledge from a large body of available datasets for power enhancement.

Another closely related strand of research focuses on the estimation of weak factors.  Influential empirical studies especially on asset pricing often give implicit yet strong evidence of weak factor models \citep{uematsu2022estimation,massacci2024instability}. {Additionally, the growing up literature on hierarchical or group factor models also supports weak factor models, as local factors load only on specific cross-sectional groups and therefore would typically be weak.}
One may refer to
\cite{ando2017clustering,choi2018multilevel,zhang2024factor,he2025factor} for empirical
evidence on such a group/hierarchical factor structure in financial and macroeconomic datasets. The definition of weak factors varies in the literature, and in this article we define  weak factors as those with  factor loadings $\bLambda^{(0)}$ scaling
sublinearly in the number $N$ of cross-section units at a slower rate, i.e., $\bLambda^{(0)\top}\bLambda^{(0)}/N^{\alpha}$  is positive definite
in the limit for some $\alpha\in (0,1)$, see also the same definitions in recent work by \cite{bai2023approximate} and \cite{choi2024high}.
The  existence of weak factors poses great challenge for large-dimensional factor modeling.  \cite{chen2024factor} proposed factor strength estimation method for factor modeling of vector and matrix-valued time series. \cite{bai2023approximate} and \cite{choi2024high} showed that the PCA estimators for weak factors are consistent under  suitable conditions on the dependence structure in the noise,  but their convergence rates would be slow  inevitably,  especially for small $\alpha$.  
\textbf{\textit{A question naturally arises: is it possible to elevate the convergence rates for the weak factors? }} In this article, we  give an affirmative answer and address this challenging problem elegantly by proposing a Transfer Principal Component Analysis (TransPCA) method, which leverages useful information from a large body of panel source datasets to enhance the estimation accuracy for weak factors. The idea originates from transfer learning in computer science, making it the last closely related strand of research.
In the past few years, the statistical theory of transfer learning in high-dimensional contexts has drawn growing attention, see for example
\cite{bastani2021predicting,li2022transfer,tian2023transfer} for supervised regressions, \cite{li2023transfer} for graphical models and \cite{li2024knowledge} for unsupervised PCA problems. 
In this work, we primarily focus on the model setup where the time horizons across the auxiliary panels and the target panel may differ, while the cross-sectional dimensions are the same for all panels. The rationale stems from the following  key considerations. Firstly,  we have rich historical time series data as source panels which would inevitably undergo  structural changes, therefore allowing the time horizons to be different provides great flexibility  
in real-world applications, see the macroeconomic data analysis in Section \ref{sec:real data} as a motivating example. Secondly, although we require the cross-sectional dimensions to be identical across target and source panels, the individual variables are allowed to be different. 
As long as the auxiliary datasets can provide useful information for estimating the weak factors of the target model, the individual variables across different panels do not need to be matched, which also provide great flexibility in real scenarios, see  the real portfolio returns data analysis in Section \ref{sec:real data} as an example.  The estimation accuracy of  weak factors for the target model would be enhanced by leveraging informative auxiliary panels  to increase the effective sample size of the target model, which is the core idea and theoretical finding of this work.

\subsection{Contributions}

The contributions of this work can be summarized as follows. Firstly, 
we propose a TransPCA method to estimate the target factor model with weak factors
by leveraging useful information from a large body of panel  datasets. In contrast to \cite{duan2024target}, our work is from a totally different perspective, focusing on the similarity directly across the loading spaces of the target and source models, which is more intuitive and flexible compared to that across the covariance matrices. In addition, our model setup is quite flexible in the sense that the individual variables across different panels need not to be completely identical. As long as the loading spaces exhibit similarity to the target one, these panels will aid in improving the estimation accuracy of the target factor model.
 Secondly, we propose a novel model setup that  the target loading space corresponding to the weak factors is shared or similar to  the loading spaces of the sources,
and integrate useful information across sources via manipulating with a novel weighted average of the  projection matrices, thereby improving the estimation accuracy of weak factors in the target model. 
Thirdly, our theoretical results reveal that when weak factors exist, the TransPCA estimators outperform the traditional PCA estimators derived solely based on the target dataset under certain assumptions, in terms of achieving remarkably superior convergence rates.  Our TransPCA method essentially improves the estimation accuracy of the target factor model by increasing the effective sample size through transfer learning, therefore TransPCA estimators  would achieve  superior or at least comparable convergence rates with those of PCA estimators in \cite{bai2003inferential} when all factors of the target model are strong. In other words, our method always achieves higher estimation accuracy than performing PCA solely on the target dataset. Furthermore, we innovatively propose a TransED method based on the maximum eigenvalue gap of the weighted average projection matrix  to determine the number of weak factors.  We also introduce two different methods for estimating the factor strengths of the target model: one  
 exclusively relying on the covariance information from the target dataset while the other  leveraging  the weighted average projection matrix.  Finally, to avoid negative transfer, we propose a criterion to select useful source panels according to a carefully designed measure of similarity level which is related to a rectified nonconvex optimization problem and provide its theoretical guarantee. {The  analysis of real macroeconomic and financial datasets illustrates the necessity of taking the availability of a large number of informative auxiliary panels into account in factor modeling and show the superiority of our proposed methods in contrast to the PCA method performed only to the target panel.}

\subsection{Organization and Notational Remarks}
The rest of this article is organized as follows. In Section \ref{sec:Methodology}, we introduce the basic framework of the transfer learning factor model and our TransPCA method, which leverages auxiliary datasets to improve the estimation accuracy of the target model. We investigate the theoretical properties of the estimated weak/strong loadings and factors in Section \ref{sec:theory}. In Section \ref{sec:deter FN}, we use Eigenvalue-Ratio method to determine target factor number and propose a Transfer Eigenvalue Difference method to consistently determine the number of weak factors. We propose a method   for estimating factor strengths and verify the consistency of the factor strength estimators in Section  \ref{sec:Estimate FS}. We further discuss the case that the informative datasets are unknown and propose a criteria and algorithm for selecting useful datasets to avoid negative transfer in Section  \ref{sec:data selection}. Numerical simulation results are reported in Section \ref{sec:simulation} and {both a real macroeconomic and a financial dataset analysis are provided in Section \ref{sec:real data}.} We conclude the paper and discuss possible future research directions in Section \ref{sec:conclusion}.

To end this section, we introduce some notations adopted throughout this paper.
For a matrix $\Ab$, we use $\Ab_{ij}$, $\Ab_{i,\cdot}$ and $\Ab_{\cdot,j}$ to denote the $(i,j)$-th entry, the $i$-th row and $j$-th column of $\Ab$, respectively. For a real symmetric matrix $\Ab$, let $\lambda_i(\Ab)$ denote the $i$-th largest eigenvalue of $\Xb$ and $\delta_i(\Ab) = \lambda_i(\Ab) - \lambda_{i+1}(\Ab)$ represent the $i$-th eigenvalue gap of $X$. Let $\mathcal{M}(\Ab)$ represent the linear space spanned by the columns of matrix $\Ab$, $\mathcal{M}(\Ab)^\perp$ denote its orthogonal complement. For matrices $\Ab$ and $\Bb$, let $\text{diag}(\Ab,\Bb)$ be a block diagonal matrices $\Ab$ and $\Bb$.
We use $\|\Ab\|_2$ and $\|\Ab\|_F$ to denote the Frobenius norm and spectral of matrix $\Ab$. Let $\Ab^\top$, $\text{tr}(\Ab)$ and $\text{rank}(\Ab)$ be the transpose, the trace and the
rank of $\Ab$, respectively. For any vector $\ba$, denote its $\ell_2$ norm as $\|\ba\|_2$. 
We use $\lfloor x \rfloor$ and $\lceil x \rceil$ to denote the largest previous and smallest following integers of $x$.
For two scalars sequences $\{a_n\}_{n\geq 1}$ and $\{b_n\}_{n\geq 1}$, we say $a_n\lesssim b_n$ ($a_n\gtrsim b_n$) if there exists a universal constant $C>0$ such that $a_n\geq C b_n$ ($a_n\leq C b_n$) and $a_n\asymp b_n$ if both $a_n\lesssim b_n$ and $a_n\gtrsim b_n$ hold. We denote the set $\{1,\ldots,T\}$ by $[T]$ for convenience. For two sets $\mathcal{A}$ and $\mathcal{B}$,  let $\mathcal{A}\setminus\mathcal{B}$ denotes the  intersection of set $\mathcal{A}$ and the complement of  $\mathcal{B}$, i.e, $\mathcal{A}\cap \mathcal{B}^{\text{c}}$. For a full column rank matrix $\Ab\in\RR^{m\times n}$, let $\Pb_{\Ab}=\Ab(\Ab^\top\Ab)^{-1}\Ab^\top$ be the orthogonal projection matrix to $\mathcal{M}(\Ab)$, $(\Pb_\Ab)^\perp=\Ib_m-\Pb_\Ab$. For convenience, denote the total sample size as 
$T = \sum_{k \in \{0\}\cup \mathcal{A}}T_k$, where $T_k$ is the sample size of the $k$-th dataset.

\section{Model setup and Methodology}\label{sec:Methodology}
In the following, we primarily introduce the factor model transfer learning setup   in Section \ref{sec:model-setup} and propose the oracle TransPCA method in Section \ref{oracle TransPCA}.
\subsection{Model setup}\label{sec:model-setup}
We observe a target panel dataset $\Xb^{(0)}=(\bX_1^{(0)},\ldots,\bX_{T_0}^{(0)})^\top$ with $T_0$ time periods and $N$ cross-sectional units. For the target panel, assume it has the following factor structure
\begin{equation}\label{VFM}
\bX_t^{(0)}= \bLambda^{(0)}\bm{f}_t^{(0)}+\bepsilon_t^{(0)},
\end{equation}
where $\bLambda^{(0)}=\left(\blambda_1^{(0)}, \dots,\blambda_{r_0}^{(0)}\right) \in \RR^{N \times r_0}$ is the factor loading matrix, $\bm{f}_t^{(0)}$ is the unobservable $r_0 \times 1$ factor  at time $t$, $\bepsilon_t^{(0)}=\left(\epsilon_{t1}^{(0)},\ldots,\epsilon_{tN}^{(0)}\right)^\top$ is $N \times 1$ idiosyncratic error and  $r_0$ is the number of factors which is assumed to be known and fixed at the moment and we will discuss the determination of $r_0$ in the following.  Factor model (\ref{VFM}) can also be written in the matrix form
\begin{equation}\label{target FM}
\Xb^{(0)}= \Fb^{(0)} {\bLambda^{(0)\top}}+\Eb^{(0)},
\end{equation}
where $\Fb^{(0)}=(\bm{f}_1^{(0)}, \bm{f}_2^{(0)}, \dots, \bm{f}_{T_0}^{(0)})^\top$, $\Eb^{(0)}=(\bepsilon_1^{(0)}, \bepsilon_2^{(0)}, \dots, \bepsilon_{T_0}^{(0)})^\top$.

Our goal is to estimate the latent factor  structure in $\Xb^{(0)}$. However, in practical applications, due to the presence of weak factors or the small time horizon $T_0$, it is challenging to achieve accurate estimation for the latent factor model in (\ref{target FM}). Let $\blambda^{(0)}_{i}$ be the $i$-th column of $\bLambda^{(0)}$, and its factor strength, characterized by $\alpha_{i}$ ranging between 0 and 1, is defined such that $$\blambda_i^{(0)\top}\blambda_i^{(0)}\asymp N^{\alpha_i}, \ \ \ i\in [r_0].$$
For a strong factor, its factor strength $\alpha_{i} = 1$, while for a weak factor, its factor strength $0<\alpha_{i} <1$. Indeed, a weak
factor may arise from two scenarios:  the factor has a weak effect on all time series or  it affects
only a subset of the series. Group factor structures would naturally fall into the latter scenario \citep{zhang2024factor,he2025factor}, where the group-specific factors affect
only a subset of time series and the loadings are
nonzero  only among specific cross-sectional groups/clusters.
Without loss of generality, we assume there exist $s$ weak factors for the target  factor model, and let
$\bLambda_w^{(0)} \in\RR^{N\times s}$ and $\bLambda_s^{(0)}\in \RR^{N\times (r_0-s)}$ be the loading matrices  corresponding to the weak factors and strong factors in $\Fb^{(0)}$, denoted 
 by $\Fb^{(0)}_w$ and $\Fb^{(0)}_s$,  respectively. 

The estimators of $\bLambda^{(0)}$ and $\Fb^{(0)}$ can be obtained by Principal Component Analysis (PCA), which is equivalent to solving   the following optimization problem
$$\underset{\bLambda^{(0)}, \ \Fb^{(0)}}{\min}\frac{1}{NT_0 }\Big\|\Xb^{(0)}-\Fb^{(0)}\bLambda^{(0)\top}\Big\|_F^2,$$
subject to $\bLambda^{(0)\top}\bLambda^{(0)}=\Db$, where $\Db$ is a diagonal matrix with diagonal elements $(\Db)_{ii}\asymp N^{\alpha_{i}}$, $i\in [r_0]$.
The solution to the above least squares problem is $$\hat{\bLambda}^{(0)}=\hat{\Qb}^{(0)}\Db^{1/2},$$ where $\hat{\Qb}^{(0)}$ is composed of the  top $r_0$ eigenvectors of the sample covariance matrix $\hat{\bSigma}^{(0)}=\frac{1}{T_0}\Xb^{(0)\top}\Xb^{(0)}$,
$\hat{\Fb}^{(0)}$ can be further obtained  by directly regressing $\Xb^{(0)}$ on $\hat{\bLambda}^{(0)}$, i.e., 
$\hat{\Fb}^{(0)}=\Xb^{(0)\top}\hat{\bLambda}^{(0)}\left(\hat{\bLambda}^{(0)\top}\hat{\bLambda}^{(0)}\right)^{-1}=\Xb^{(0)\top}\hat{\bLambda}^{(0)}\Db^{-1}.$
We first assume that $\Db$ is given in advance, and we will discuss the estimation of  $\Db$ in Section \ref{sec:Estimate FS}. 

Note that both $\Fb^{(0)}$ and $\bLambda^{(0)}$ are unobservable/unidentifiable and can be replaced by $\Fb^{(0)} (\Hb^{-1})^\top$ and $\bLambda^{(0)}\Hb$ for any invertible matrix $\Hb\in \RR^{r_0\times r_0}$ without changing the model. Although the loading matrix $\bLambda^{(0)}$ is not identifiable, the space spanned by the columns of $\bLambda^{(0)}$, called loading space and denoted by $\mathcal{M}(\bLambda^{(0)})$, is uniquely defined. Thus, one of the primary goals in this study is to estimate the target loading (factor) space rather than the loading (factor) matrix.

The existence of weak factors may result in slower convergence rates of the PCA estimators 
compared with those derived under the pervasive/strong factor condition in the literature \citep{fan2013large}. To improve the estimation accuracy of the target dataset, we aim to leverage
additional information from auxiliary panel datasets $\Xb^{(k)}$, $k\in \mathcal{A} \subseteq [K]$. 
Consider an oracle case that  all informative datasets $\Xb^{(k)}$ are given, i.e., $\mathcal{A}$ is known in advance and
each informative dataset also has an approximate latent factor structure
\begin{equation*}\label{source FM}
\Xb^{(k)}= \Fb^{(k)} {\bLambda^{(k)\top}}+\Eb^{(k)}, \quad k\in \mathcal{A}\ ,
\end{equation*}
where $\Xb^{(k)} \in \RR^{T_k \times N}$, $\Fb^{(k)} \in \RR^{T_k \times r_k}$ is the factor matrix, $\bLambda^{(k)} \in \RR^{N \times r_k}$ is the factor loading matrix and $\Eb^{(k)}\in \RR^{T_k \times N}$ is  the idiosyncratic error matrix,  $r_k$ is the number of factors which is assumed to be known and fixed and its estimation will also be discussed in the following sections. 
For any $k \in \mathcal{A}$, 
let $\bLambda^{(k)}_\mathcal{T}\in\RR^{N\times s}$ be the matrix composed of $s$ columns of 
$\bLambda^{(k)}$  whose spanned space is informative for estimating $\bLambda^{(0)}_w$, while  $\bLambda^{(k)}_{\mathcal{T}^c}\in\RR^{N\times (r_k-s)}$
be the matrix composed of the remaining columns of $\bLambda^{(k)}$. 
We can perform traditional PCA method on $\Xb^{(k)}$ to obtain estimators of the loadings $\bLambda^{(k)}$ and factors $\Fb^{(k)}$.
Recall that our primary goal  is to  estimate the target loading (factor) space. For the auxiliary datasets, we are particularly interested in how to leverage useful information such that their loading (factor) spaces contribute to improving the estimation of loading (factor) space of the target dataset. To evaluate the informative level of $\mathcal{M}(\bLambda_{\mathcal{T}}^{(k)})$, we assume for any $k\in \mathcal{A}$, for some $\varepsilon>0$,
\begin{equation}\label{similarity}
\Big\|\Pb_{\bLambda_\mathcal{T}^{(k)}}-\Pb_{\bLambda_w^{(0)}}\Big\|_F\leq \varepsilon.
\end{equation}

\begin{remark}
For convenience, we assume that the loading spaces of informative auxiliary datasets are similar to  the loading space corresponding to weak factors of the target model.
It is likely that the informative auxiliary datasets are also beneficial for estimating the loading space corresponding to  strong factors of the target model. We point out that our proposed method can be easily extended to deal with this case and this assumption here is only for simplicity of illustration.
\end{remark}
\subsection{Oracle Transfer Principal Component Analysis Procedure}\label{oracle TransPCA}
In this section, we introduce our Transfer Principal Component Analysis (TransPCA) method  to estimate the target factor model by leveraging useful information from a large body of panel  datasets. We first consider the oracle case that all the informative datasets are given, i.e, the index set $\cA$ is known in advance. To distinguish it from the procedure equipped with useful dataset selection capability in Section \ref{sec:data selection}, we refer to the TransPCA method in this section as Oracle TransPCA.
We utilize a weighted average projection matrix to aggregate useful information from auxiliary datasets, which is the primary idea and contribution of the current work. Specifically, define
\begin{equation}\label{sigma-w}
\begin{aligned}
\Pb^{w}&=\frac{1}{T}\sum_{k\in\{0\}\cup\mathcal{A}} T_k \Pb_{\bLambda^{(k)}}=\frac{1}{T}\left(T_0\left(\Pb_{\bLambda_w^{(0)}}+\Pb_{\bLambda_s^{(0)}}\right)+\sum_{k\in  \mathcal{A} }T_k\left(\Pb_{\bLambda_{\mathcal{T}}^{(k)}} +\Pb_{\bLambda_{\mathcal{T}^c}^{(k)}} \right)
\right)\\
&=\Pb_{\bLambda_{w}^{(0)}}+\underbrace{\frac{1}{T}\sum_{k\in\mathcal{A}}T_k\left(\Pb_{\bLambda_{\mathcal{T}}^{(k)}}-\Pb_{\bLambda_w^{(0)}}\right)}_\text{Similarity Term}+\underbrace{\frac{1}{T}\left(T_0 \Pb_{\bLambda_s^{(0)}}+\sum_{k\in\mathcal{A}}T_k \Pb_{\bLambda_{\mathcal{T}^c}^{(k)}}\right)}_\text{Difference Term},
\end{aligned}
\end{equation}
where $T = \sum_{k \in \{0\}\cup \mathcal{A}}T_k$. The ``Similarity Term'' in (\ref{sigma-w}) primarily measures the similarity between the target dataset and the auxiliary datasets, which determines the degree that how transferable and informative are the auxiliary datasets. The ``Difference Term'' consists of the ``private" space terms of the target  and the auxiliary datasets which could be quite different. From equation (\ref{sigma-w}), we can conclude that if the ``Similarity Term'' is sufficiently small and $\Pb_{\bLambda_w^{(0)}}$ is the dominant term compared with the ``Difference Term'', the eigenvectors of $\Pb^w$ can be sufficiently close to  $\bLambda_w^{(0)}$. 

Inspired by equation (\ref{sigma-w}), 
to leverage useful information from the auxiliary datasets, 
 we introduce the following two-step Oracle TransPCA procedure to estimate $\bLambda^{(0)}$ and $\Fb^{(0)}$:
\\
\textbf{The First Step:} for each
$\ k \in \{0\} \cup \mathcal{A}$, 
let $\hat{\Qb}^{(k)}$  be the leading  $r_k$ eigenvectors of the sample covariance matrix   $$\hat{\bSigma}^{(k)}=\frac{1}{NT_k}\Xb^{(k)\top}\Xb^{(k)},$$
and let $\hat{\Qb}_w^{(0)}$ be the leading $s$ eigenvectors of the weighted average projection matrix  defined as 
\begin{equation*}
\hat{\Pb}^w = \frac{1}{T}\sum_{k \in \{0\}\cup \mathcal{A}}T_k \cdot\hat{\Qb}^{(k)} \hat{\Qb}^{(k)\top}.
\end{equation*}
Next, we obtain the estimator of the weak loading matrix, defined as  $\hat{\bLambda}^{(0)}_w=\hat{\Qb}_w^{(0)}\hat{\Db}_w^{1/2}$,  where  $\hat{\Db}_w$ is a diagonal matrix with diagonal elements reflecting the strength of weak factors, as will be further discussed in Section \ref{sec:Estimate FS}. Then, we estimate the weak factors by $\hat{\Fb}_w^{(0)}=\Xb^{(0)}\hat{\bLambda}^{(0)}_w\hat{\Db}_w^{-1}$, which is obtained by regressing $\Xb^{(0)}$ on $\hat{\bLambda}^{(0)}_w$.

{After obtaining $\hat{\bLambda}_w^{(0)}$, we then project $\hat{\bLambda}_w^{(0)}$ out from $\Xb^{(0)}$, i.e,   we set $\hat{\Yb}^{(0)}=\Xb^{(0)}\left(\Ib_N-\hat{\Qb}_w^{(0)}\hat{\Qb}_w^{(0)\top}\right)$ and then perform PCA on $\hat{\Yb}^{(0)}$ to estimate $\bLambda_s^{(0)}$. }
\\
\textbf{The Second Step:} we acquire $\hat{\Qb}_s^{(0)}$ by taking the leading $(r_0-s)$ eigenvectors of the following projected sample covariance matrix
\begin{equation}\label{projected cov}
\frac{1}{NT_0}\hat{\Yb}^{(0)\top}\hat{\Yb}^{(0)}=\left(\Ib_N-\hat{\Qb}_w^{(0)}\hat{\Qb}_w^{(0)\top}\right)\hat{\bSigma}^{(0)}\left(\Ib_N-\hat{\Qb}_w^{(0)}\hat{\Qb}_w^{(0)\top}\right),
\end{equation}
and estimate the loadings corresponding to strong factors by 
$\hat{\bLambda}^{(0)}_s=\sqrt{N}\cdot\hat{\Qb}_s^{(0)}$. Similarly, the estimator of  strong factors $\hat{\Fb}^{(0)}_s=\frac{1}{N}\left(\Xb^{(0)}-\hat{\Fb}_w^{(0)}\hat{\bLambda}^{(0)\top}_w\right)\hat{\bLambda}^{(0)}_s$ is obtained  by directly regressing $(\Xb^{(0)}-\hat{\Fb}_w^{(0)}\hat{\bLambda}^{(0)\top}_w)$ on $\hat{\bLambda}_s^{(0)}$. At last, we summarize the  Oracle TransPCA procedure in Algorithm \ref{algorithm:1} for better illustration.

\begin{algorithm}
\caption{Oracle Transfer Principal Component Analysis Procedure}\label{algorithm:1}
\begin{algorithmic}[1]
\REQUIRE Datasets $\Xb^{(k)}$ and the corresponding factor numbers $r_k$ for $k\in \{0\}\cup \mathcal{A}$; weak factor strength $\Db_w$; the number of weak factors  $s$. 
\STATE Perform PCA on each dataset $\big\{\Xb^{(k)},\ k\in \{0\}\cup \mathcal{A}\big\}$, to obtain $\hat{\Qb}^{(k)}$.
\STATE Define $\hat{\Qb}^{(0)}_w$ as the leading $s$ eigenvectors of $\hat{\Pb}^w = \frac{1}{T}\sum_{k \in \{0\}\cup \mathcal{A}} T_k \Pb_{\hat{\Qb}^{(k)}}$, and   let $\hat{\bLambda}^{(0)}_w = \hat{\Qb}_w^{(0)}\Db^{1/2}_w$.
\STATE  Obtain $\hat{\Fb}_w^{(0)}=\Xb^{(0)}\hat{\bLambda}^{(0)}_w\hat{\Db}_w^{-1}$ by regressing $\Xb^{(0)}$ on $\hat{\bLambda}_w^{(0)}$.
\STATE Define $\hat{\Qb}^{(0)}_s$ as the leading $(r_0-s)$ eigenvectors of the projected
sample covariance matrix  in (\ref{projected cov}), and let $\hat{\bLambda}^{(0)}_s = \sqrt{N} \cdot \hat{\Qb}_s^{(0)}$.
\STATE Obtain $\hat{\Fb}^{(0)}_s=\frac{1}{N}\left(\Xb^{(0)}-\hat{\Fb}_w^{(0)}\hat{\bLambda}^{(0)\top}_w\right)\hat{\bLambda}^{(0)}_s$ by regressing $\Xb^{(0)}-\hat{\Fb}_w^{(0)}\hat{\bLambda}_w^{(0)\top}$ on $\hat{\bLambda}_s^{(0)}$.

\ENSURE The estimators of loadings and factors for the target model, $\hat{\bLambda}_w^{(0)}$, $\hat{\bLambda}_s^{(0)}$, $\hat{\Fb}_w^{(0)}$ and $\hat{\Fb}_s^{(0)}$.
\end{algorithmic}
\end{algorithm}



\section{Theoretical Results} \label{sec:theory}
In this section, we investigate the theoretical properties of the proposed TransPCA estimators. We first give some assumptions.

\begin{assumption}\label{factors}
(Factors) For any $k\in \{0\}\cup [K]$, assume that
\noindent 
\begin{enumerate}
    \renewcommand{\labelenumi}{(\theenumi)}
    \item 
    $E\|\bm{f}_t^{(k)}\|^4 \leq M <\infty$;
    \item $\frac{\Fb^{(k)\top}\Fb^{(k)}}{T_k}\overset{P}{\longrightarrow}\bSigma_{F}^{(k)}$ as $T_k\rightarrow \infty$, where $\bSigma_{F}^{(k)}$ is a $r_k\times r_k$ positive definite matrix with bounded distinctive eigenvalues.   
\end{enumerate}
\end{assumption}

\begin{assumption}\label{loadings}
(Factor loadings)

\begin{enumerate}
    \renewcommand{\labelenumi}{(\theenumi)}
    \item $\|\bLambda^{(k)}\|_{\text{max}} \leq M< \infty$, $\forall k\in \{0\}\cup [K]$;
    \item 
for the target dataset, $\Bigg\|\Db^{-1/2}\bLambda^{(0)\top}\bLambda^{(0)}\Db^{-1/2}-\Ib_{r_0}\Big\| \rightarrow 0$, where $\Db=\text{diag}(\Db_s,\Db_w)$ is a $r_0\times r_0$ diagonal matrix, $(\Db_s)_{ii}\asymp N$ for $i \in[r_0-s]$, $(\Db_w)_{ii}\asymp N^{\alpha_i}$ for $i\in[s]$, and $0<\alpha_s\leq\cdots\leq\alpha_1< 1$;


    \item for the auxiliary datasets, $\big\|N^{-1}\bLambda^{(k)\top}\bLambda^{(k)}-\Ib_{r_k}\big\|\rightarrow 0$, $\forall k\in[K]$.
\end{enumerate}
\end{assumption}

\begin{assumption}\label{errors}
(Idiosyncratic errors)
    For any $k\in \{0\}\cup [K]$,     $\big\|\Eb^{(k)}\big\|^2=O_p\left(N \vee T_k\right)$.

\end{assumption}

\begin{assumption}\label{corr}
(Weak time and cross-section correlations)
    \renewcommand{\labelenumi}{(\theenumi)}
    \noindent 
\begin{enumerate}
    \item 
    For each $t$, $E\Big\|\Db^{-1/2}\sum_{i}\bLambda^{(0)}_{i,\cdot}\epsilon^{(0)}_{ti}\Big\|^2 \leq M$, $E\Big\|N^{-1/2}\sum_i\bLambda_{i,\cdot}^{(k)}\epsilon_{ti}^{(k)}\Big\|^2\leq M$, $\forall k\in[K]$;
    \item 
    for each $i$, $E\Big\|T_k^{-1/2}\sum_{t}\bm{f}_t^{(k)}\epsilon^{(k)}_{ti}\Big\|^{2} \leq M$,  $\forall k\in \{0\}\cup [K]$.
     \end{enumerate}
\end{assumption}

Assumption \ref{factors} and Assumption \ref{loadings}.(1) are standard for factor models, see for example \cite{bai2002determining} and \cite{bai2003inferential}. Assumption \ref{loadings}.(2) allows for the existence of weak factors by assuming the limit of $\Db^{-1/2}\bLambda^{(0)\top}\bLambda^{(0)}\Db^{-1/2}$ to be the identity matrix. Assumption \ref{loadings}.(2) requires that the factors of all auxiliary datasets are strong, which is mainly for simplifying the theoretical analysis and can be further relaxed.  A similar assumption can also be found in 
\cite{duan2024target}.
Assumption \ref{errors} imposes an upper bound on the maximum eigenvalues of $\Eb^{(k)\top}\Eb^{(k)}$, ruling out the possibility that the error matrix $\Eb^{(k)}$ contains common factors. For $\Eb^{(k)}$ whose elements are
$i.i.d.$  random variables with uniformly bounded fourth moments, Assumption \ref{errors} is easily satisfied. Assumption \ref{corr} allows weak cross-sectional dependence in the idiosyncratic components and 
weak dependence between factors and idiosyncratic errors. Assumptions \ref{factors}-\ref{corr} are  similar to those in \cite{bai2002determining}, \cite{bai2003inferential}, \cite{ahn2013eigenvalue}, and 
the same as those in \cite{bai2023approximate} which also account for weak factors. To derive the theoretical properties of our TransPCA estimators, we also need  the following assumptions.

\begin{assumption}\label{Identifiability}
Given the decomposition $\Pb_{\bLambda^{(k)}} = \Pb_{\bLambda_\mathcal{T}^{(k)}} + \Pb_{\bLambda_{\mathcal{T}^c}^{(k)}}$ for $k \in \mathcal{A}$, for the informative and transferable loading spaces, we assume that $\Big\|\Pb_{\bLambda_\mathcal{T}^{(k)}}-\Pb_{\bLambda_w^{(0)}}\Big\|_F\leq \varepsilon$. As for the ``Difference Term'' in (\ref{sigma-w}), we assume that there exists some  constant $\kappa>0$ such that 
\begin{equation}\label{equ:non-trans-gap}
\frac{1}{T}\left\|T_0\Pb_{\bLambda_s^{(0)}}+
\sum_{k \in  \mathcal{A}} T_k\Pb_{\bLambda_{\mathcal{T}^c}^{(k)}}\right\| \leq 1-\kappa, 
\end{equation}
where $T = \sum_{k \in \{0\}\cup \mathcal{A}} T_k$.
\end{assumption}

\begin{assumption}\label{NT-relation}
As $N$, $T\rightarrow \infty$,
(1) $\frac{N}{\sqrt{T_0}N^{\alpha_s}}\rightarrow 0$; (2) $\frac{T_0\sqrt{N}}{T\sqrt{N^{\alpha_s}}}\rightarrow 0$.
\end{assumption}

Assumption \ref{Identifiability}
is crucial for distinguishing transferable subspaces from non-transferable ones.
It requires that the non-transferable projection matrices $\Pb_{\bLambda_{\mathcal{T}^c}^{(k)}}$ do not jointly
contribute to some large eigenvalues of $\Pb^w$ in (\ref{sigma-w}), 
and a similar assumption can also be found
 in \cite{duan2024target}.
Assumption \ref{Identifiability} 
 ensures that $\Pb_{\bLambda_w^{(0)}}$ is the dominant term in (\ref{sigma-w}) when $\varepsilon \rightarrow 0$. 
Assumption \ref{NT-relation} is a scaling condition on $N$, $T_0$ and $T$, while Assumption \ref{NT-relation}.(1) is more relaxed than that in \cite{bai2023approximate} and \cite{duan2024target} which assume that ${N}/{(T_0 N^{\alpha_s})}\rightarrow 0$.

We first assume that the number of weak factors in target model $s$, the number of factors $r_k, k\in \{0\}\cup \mathcal{A }$ and the factor strengths $\alpha_i, i\in [s]$ for the target model are known in advance. In Section \ref{sec:deter FN} and Section \ref{sec:Estimate FS}, we will provide 
consistent estimators for factors numbers and factor strengths, respectively.
Theorem \ref{th1} shows the  error rates in terms of estimating the loading and factor spaces corresponding to  weak factors.

\begin{theorem}(General stochastic order for estimation errors)\label{th1}
Under Assumptions \ref{factors} to \ref{Identifiability} and  \ref{NT-relation}.(1) with the number of factors $s$, $r_k, k\in \{0\}\cup \mathcal{A }$ fixed and given, as $N$, $ \underset{k\in\{0\}\cup \mathcal{A}}{\min}T_k \rightarrow \infty$,  we have
$$\frac{1}{\sqrt{N}}\bigg\|\hat{\bLambda}_w^{(0)}-\bLambda_w^{(0)}\Hb_w\bigg\|=\frac{\sqrt{N^{\alpha_1}}}{\sqrt{N}}\cdot O_p\left(\frac{\sqrt{T_0}N}{TN^{\alpha_s}}+\frac{T_0}{TN^{\alpha_s}}+\sum_{k\in\mathcal{A}}\left(\frac{\sqrt{T_k}}{T}+\frac{T_k}{TN}\right)+\frac{1}{T}\sum_{k\in \mathcal{A}}T_k \varepsilon \right),$$
$$\frac{1}{\sqrt{T_0}}\bigg\|\hat{\Fb}_w^{(0)}-\Fb_w^{(0)}\left(\Hb_w^{-1}\right)^{\top}\bigg\|=\frac{\sqrt{N}}{\sqrt{N^{\alpha_s}}}\cdot O_p\left(\frac{\sqrt{T_0}N}{TN^{\alpha_s}}+\frac{T_0}{TN^{\alpha_s}}+\sum_{k\in\mathcal{A}}\left(\frac{\sqrt{T_k}}{T}+\frac{T_k}{TN}\right)+\frac{1}{T}\sum_{k\in \mathcal{A}}T_k \varepsilon\right)+O_p\left(\frac{1}{\sqrt{N^{\alpha_s}}}\right).$$
\end{theorem}
{Theorem \ref{th1} presents general stochastic order of the errors in terms of estimating the weak loadings and factors spaces. By imposing certain condition on the similarity level $\varepsilon$, as well as additional scaling condition on $N$, $T_0$ and $T$, we would expect faster convergence rate compared with those derived in the existing literature.}

\begin{corollary}\label{cor:1}
Under Assumptions \ref{factors} to  \ref{NT-relation} with the number of factors $s$, $r_k, k\in \{0\}\cup \mathcal{A }$ fixed and given, further assume that 
$\varepsilon=o\left({1}/{\sqrt{N^{1-\alpha_s}}}\right)$, then  we have
as $N$, $ \underset{k\in\{0\}\cup \mathcal{A}}{\min}T_k \rightarrow \infty$, 
$$\frac{1}{\sqrt{N}}\bigg\|\hat{\bLambda}_w^{(0)}-\bLambda_w^{(0)}\Hb_w\bigg\|=o_p(1), \quad \frac{1}{\sqrt{T_0}}\bigg\|\hat{\Fb}_w^{(0)}-\Fb_w^{(0)}\left(\Hb_w^{-1}\right)^{\top}\bigg\|=o_p(1).$$
\end{corollary}
In corollary \ref{cor:1}, with additional  Assumption \ref{NT-relation}.(2) and $\varepsilon=o\left({1}/{\sqrt{N^{1-\alpha_s}}}\right)$,  we  obtain consistent estimators for the weak loadings and factors up to a rotation. 
\begin{remark}
When all  factors are strong for the target model, i.e., $\alpha_i=1$ for $i\in[s]$, Theorem \ref{th1} yields
$$\frac{1}{\sqrt{N}}\bigg\|\hat{\bLambda}_w^{(0)}-\bLambda_w^{(0)}\Hb_w\bigg\|=O_p\left(\sum_{k\in\{0\}\cup \mathcal{A}}\left(\frac{\sqrt{T_k}}{T}+\frac{T_k}{TN}\right)+\frac{1}{T}\sum_{k\in \mathcal{A}}T_k \varepsilon \right).$$
If we assume that $\varepsilon=o\left({1}/{\sqrt{T}}+{1}/{N}\right)$, the convergence rate of the estimated loading matrix is $\min\{\sqrt{T},N\}$.
When  $T_0\asymp T$, this rate is the same as the convergence rate derived in \cite{bai2003inferential} and \cite{bai2023approximate} under the strong factors assumption, while if $T_0\ll T$, and $\sqrt{T}\asymp N$ or $\sqrt{T}\ll N$, our TransPCA estimator would achieve a faster convergence rate compared with that in \cite{bai2003inferential} and \cite{bai2023approximate}.
\end{remark}

The following proposition  presents the convergence rate of our TransPCA estimator for the target loading space when all factors are weak. In this case, the second step of our TransPCA procedure  can indeed be omitted.

\begin{proposition}\label{pro:1}
Under Assumptions \ref{factors} to  \ref{NT-relation} with the number of factors $s$ and  $r_k$ for $k\in \mathcal{A }$ fixed and given, further assume that the factor strength matrix $\Db=\Db_w$, i.e, all factors are weak, then as $N$, $ \underset{k\in\{0\}\cup \mathcal{A}}{\min}T_k \rightarrow \infty$, we have
$$\frac{1}{\sqrt{N}}\bigg\|\hat{\bLambda}_w^{(0)}-\bLambda_w^{(0)}\Hb_w\bigg\|=\frac{\sqrt{N^{\alpha_1}}}{\sqrt{N}}\cdot O_p\left(\frac{\sqrt{T_0}\sqrt{N^{1+\alpha_1}}}{TN^{\alpha_s}}+\frac{T_0}{TN^{\alpha_s}}+\frac{N}{TN^{\alpha_s}}+\sum_{k\in\mathcal{A}}\left(\frac{\sqrt{T_k}}{T}+\frac{T_k}{TN}\right)+\frac{1}{T}\sum_{k\in \mathcal{A}}T_k \varepsilon \right).$$
\end{proposition}

\begin{remark}
When all factors are weak for the target model, i.e., $\Db=\Db_w$, \cite{bai2023approximate} shows that 
\begin{equation}\label{bai2023}
\frac{1}{\sqrt{N}}\bigg\|\hat{\bLambda}_w^{(0)}-\bLambda_w^{(0)}\Hb_w\bigg\|=O_p\left(\frac{\sqrt{N^{\alpha_1}}}{\sqrt{T_0}\sqrt{N^{\alpha_s}}}+\frac{1}{\sqrt{N^{1+\alpha_s}}}\right).
\end{equation}
We further consider the case where the target sample size is much smaller than the whole sample size, i.e., $T_0\ll T$, which is typically the case for real financial/economic application as we have rich historical data as auxiliary panels. Proposition  \ref{pro:1} yields that
\begin{equation}\label{trans-loading}
\frac{1}{\sqrt{N}}\bigg\|\hat{\bLambda}_w^{(0)}-\bLambda_w^{(0)}\Hb_w\bigg\|=\frac{\sqrt{N^{\alpha_1}}}{\sqrt{N}}\cdot O_p\left(\frac{\sqrt{T_0}\sqrt{N^{1+\alpha_1}}}{TN^{\alpha_s}}+\frac{T_0}{TN^{\alpha_s}}+\frac{1}{\sqrt{T}}+\frac{1}{N}+\varepsilon\right).
\end{equation}
By comparing (\ref{bai2023}) and (\ref{trans-loading}),  when
$T_0 \ll N$ or $T_0\asymp N$, $\varepsilon=o\left(\max\Big\{\frac{\sqrt{T_0}\sqrt{N^{1+\alpha_1}}}{TN^{\alpha_s}},\frac{T_0}{TN^{\alpha_s}},\frac{1}{\sqrt{T}},\frac{1}{N}\Big\}\right)$ and Assumption \ref{NT-relation} holds, our TransPCA estimator would achieve a faster convergence rate than that of the estimator 
in \cite{bai2023approximate}, which is obtained by performing PCA directly on the target dataset.
If further assume that $$\frac{\sqrt{T_0}\sqrt{N^{1+\alpha_1}}}{TN^{\alpha_s}}+\frac{T_0}{TN^{\alpha_s}}=o\left(\frac{1}{\sqrt{T}}+\frac{1}{N}\right),$$ then we would achieve a faster convergence rate than that derived in \cite{bai2003inferential},
which performs PCA to the target dataset $\Xb^{(0)}$ presuming that all factors are strong.

\end{remark}

From Theorem \ref{th1} and the above discussion, we can conclude that as long as there exists a useful auxiliary dataset $\Xb^{(k)}$
with  sample size $T_k\gg T_0$, we would effectively extract useful information from the auxiliary datasets by utilizing the weighted projection matrix $\hat{\Pb}^w$, thereby significantly improving the estimation accuracy. Essentially, our TransPCA method improves estimation accuracy by increasing the effective sample size, which is quite straightforward and reasonable. In the following we establish the convergence rates of the estimated  loadings and factors corresponding to strong factors for our TransPCA method and
for ease of notation, we denote $$\delta=\frac{\sqrt{N^{\alpha_1}}}{\sqrt{N}}\left(\frac{\sqrt{T_0}N}{TN^{\alpha_s}}+\frac{T_0}{TN^{\alpha_s}}+\frac{1}{\sqrt{T}}+\frac{1}{N}+\varepsilon\right).$$

\begin{theorem}\label{th2}
Under Assumptions \ref{factors} to \ref{NT-relation} with the number of factors $s$, $r_k, k\in \{0\}\cup \mathcal{A }$ fixed and given,  further assume that $\varepsilon=o\left({1}/{\sqrt{N^{1-\alpha_s}}}\right)$, then as $N$, $ \underset{k\in\{0\}\cup \mathcal{A}}{\min}T_k \rightarrow \infty$,  we have

$$\frac{1}{\sqrt{N}}\bigg\|\hat{\bLambda}_s^{(0)}-\bLambda_s^{(0)}\Hb_s\bigg\|=O_p\left(\frac{\sqrt{N}}{\sqrt{N^{\alpha_s}}}\cdot\delta+\frac{1}{\sqrt{T_0}}+\frac{\sqrt{N^{\alpha_1}}}{\sqrt{N^{1+\alpha_s}}}\right),$$
$$\frac{1}{\sqrt{T_0}}\bigg\|\hat{\Fb}_s^{(0)}-\Fb_s^{(0)}\left(\Hb_s^{\top}\right)^{-1}\bigg\|=O_p\left(\frac{\sqrt{N}}{\sqrt{N^{\alpha_s}}}\cdot\delta+\frac{1}{\sqrt{T_0}}+\frac{\sqrt{N^{\alpha_1}}}{\sqrt{N^{1+\alpha_s}}}\right).$$

\end{theorem}

Theorem \ref{th2} 
establishes the convergence rate of the estimated strong loadings and factors up to rotation. When all factors are strong in the target model, the convergence rates of the estimators for strong factors and loadings  are both $\min\{\sqrt{T_0},\sqrt{N} \}$ which are exactly the same with those derived in \cite{bai2002determining}. As we assume that the loading spaces of informative auxiliary datasets are only similar
to the loading space corresponding to weak factors of the target dataset, we can at most achieve  the convergence rate of the form $\sqrt{T_0}$ rather than $\sqrt{T}$. When $\sqrt{T_0}\ll N$, under some other reasonable/mild  assumptions, our TransPCA estimators can achieve  faster convergence rates than those in \cite{bai2023approximate}. Indeed,
Proposition 6 in \cite{bai2023approximate} established 
the convergence rate $\Big(\sqrt{{N^{1-\alpha_s}}/{T_0}}+{1}/{\sqrt{N^{1+\alpha_s}}}\Big)$ for estimating the loading matrix $\bLambda^{(0)}$. When 
$\sqrt{T_0}\ll N$, under additional assumptions on the scaling among $N$, $T$ and $T_0$, for instance, assume that $$\max\left\{\frac{\sqrt{T_0}N}{\sqrt{T}N^{\alpha_s}}, \frac{\sqrt{TN^{\alpha_s}}}{\sqrt{T_0 N}}, \frac{\sqrt{TN^{\alpha_1}}}{N}  \right\}
 \rightarrow 0,$$ then the convergence rate of our TransPCA estimator $\hat{\bLambda}_s^{(0)}$ is $\sqrt{{N^{1-\alpha_s}}/{T}}$, which is faster than the convergence rate of the PCA estimator, $\sqrt{{N^{1-\alpha_s}}/{T_0}}$.

\section{Determining the number of factors}\label{sec:deter FN}
Thus far we have assumed that the number of weak  factors in the target model $s$ and the number of factors $r_k$, $k\in \{0\}\cup \mathcal{A}$, are known. In practice, the number of factors is unknown and needs to be estimated. In this section, we  discuss on determining the factor numbers $r_k$ and $s$. 

\subsection{Estimation of the number of factors}

There is a large amount of literature investigating eigenvalue-based
methods for selecting the number of  factors in approximate factor model.
For auxiliary dataset $\big\{\Xb^{(k)}, k\in\mathcal{A}\big\}$, the number of factors $r_k$ can be estimated by IC or PC criteria in \cite{bai2002determining}, Eigenvalue Ratio (ER) method in \citet{lam2012factor,ahn2013eigenvalue}, eigenvalue difference
criterion in \cite{onatski2010determining}. Under the strong factor assumption, we can obtain consistent estimator of $r_k$ by either method mentioned above.

For the target dataset $\Xb^{(0)}$, due to the presence of weak factors, we need an additional assumption to ensure $r_0$ can be consistently estimated. 

\begin{assumption}\label{target-eigengap-noise}
For some positive and finite real number $C$ and some $c\in (0,1]$, assume  that
$$\lambda_{\lfloor c(N\land T_0)\rfloor }\left(\Eb^{(0)\top}\Eb^{(0)}/(N\lor T_0)\right)\geq C+o_p(1).$$
\end{assumption}
Assumption \ref{target-eigengap-noise} is a sufficient condition for our theoretical analysis,  indicating that the first  largest $\lfloor c(N\land T_0)\rfloor$  eigenvalues of $\Eb^{(0)\top}\Eb^{(0)}/(N\lor T_0)$ are bounded away from zero.
We use the Eigenvalue-Ratio (ER) method introduced by \cite{ahn2013eigenvalue} to estimate $r_0$. In detail, for a prespecified  maximum number of factors $r_{\max}$, we estimate $r_0$ by
$$\hat{r}_0=\underset{1\leq r\leq r_{\max}}{\arg\max}\frac{\lambda_r\left(\hat{\bSigma}^{(0)}\right)}{\lambda_{r+1}\left(\hat{\bSigma}^{(0)}\right)}.$$

Theorem \ref{th3} establishes the consistency of the ER estimator of the number of factors  in the target model.
\begin{theorem}\label{th3}
Under Assumptions \ref{factors} to \ref{corr}, \ref{NT-relation}.(1), \ref{target-eigengap-noise} and the  condition $\alpha_s>1/2$,  with $r_0\geq 1$, 
we have $P\left(\hat{r}_0=r_0\right)\rightarrow 1$ as $\min\{N,T_0\}\rightarrow \infty$, for any $r_{\max}\in(r_0,\lfloor c(N\land T_0)\rfloor -r_0-1]$.
\end{theorem}
The condition $\alpha_s>1/2$ requires that the  weakest factors are not excessively weak, ensuring a sufficient gap between the signal and the noise parts. The same condition also appears in \cite{uematsu2022estimation} and  Assumption $C'$ of \cite{bai2023approximate}.

\subsection{Estimation of the number of weak factors of the target model}
We now introduce a Transfer Eigenvalue Difference (TransED) method for determining the number of weak factors based on the maximum eigengap, similar as that in  \cite{onatski2010determining}. The following theorem illustrates the property of the eigenvalues of the weighted average projection matrix $\Pb^w$ in Section \ref{sec:Methodology}.

\begin{theorem}\label{lem:projection-eigen}
Under Assumption \ref{Identifiability}, there exists a positive constant $C$ independent of $N$ and $T$ such that
$$\lambda_j\left(\Pb^w\right)=1+C\cdot\varepsilon, \quad j\in[s], $$
$$\lambda_{s+1}\left(\Pb^w\right)\leq 1-\kappa+C\cdot\varepsilon.$$
\end{theorem}

As shown in Theorem  \ref{lem:projection-eigen}, there is a gap between the $s$-th eigenvalue and the $(s+1)$-th eigenvalue of $\Pb^w$. With this in mind, our TransED method primarily relies on  the difference between two adjacent eigenvalues of $\hat{\Pb}^w$, to determine the number of  weak factors $s$. For a prespecified  maximum number of factors $s_{\max}$, we define the estimator of the number of weak factors $s$ as
$$\hat{s}=\underset{1\leq j \leq s_{\max}}{\arg\max}\left(\lambda_j\left(\hat{\Pb}^w\right)-\lambda_{j+1}\left(\hat{\Pb}^w\right)\right).$$
 The following theorem establishes the consistency of the  estimator $\hat{s}$ for the weak factor number.
\begin{theorem}\label{th4}
Under Assumptions  \ref{factors} to  \ref{NT-relation}, further assume that  $\varepsilon=o\left({1}/{\sqrt{N^{1-\alpha_s}}}\right)$, and there exists a constant
 $c\in(0,1]$ such that
$\underset{s< j \leq \lfloor cN \rfloor }{\max}\left(\lambda_j\left(\Pb^w\right)-\lambda_{j+1}\left(\Pb^w\right)\right)<\kappa$. As $N$, $ \underset{k\in\{0\}\cup \mathcal{A}}{\min}T_k \rightarrow \infty$,  we have
$P\left(\hat{s}=s\right)\rightarrow 1$.

\end{theorem}

\section{Estimation of the factor strengths in target model}\label{sec:Estimate FS} 

 In this section, we introduce two methods for estimating the factor strengths for the target factor model: one based solely on the target dataset and the other associated with the TransPCA procedure, borrowing idea from \cite{chen2024factor}. 
 As mentioned earlier, the factor loading $\bLambda^{(0)}$ can be written as $\bLambda^{(0)} = \Qb^{(0)}\Db^{1/2}$, 
and $\Db=\text{diag}(\Db_s,\Db_w)$ is a diagonal matrix, where $(\Db_s)_{ii}\asymp N$, $i\in[r_0-s]$, $(\Db_w)_{ii}\asymp N^{\alpha_i}$, $i\in[s]$.

 Define $\hat{\Mb} = \hat{\Qb}^{(0)\top} \left(\frac{1}{T_0}\Xb^{(0)\top}\Xb^{(0)}\right)\hat{\Qb}^{(0)}$, where
$\hat{\Qb}^{(0)}$ is a column orthogonal matrix composed of the leading $r_0$ eigenvectors of $\Xb^{(0)\top} \Xb^{(0)}/{T_0}$. Then, we have
 \begin{equation*}
 \begin{aligned}
 \hat{\Mb}&=\hat{\Qb}^{(0)\top}\Qb^{(0)}\Db^{1/2}\left(\frac{1}{T_0}\Fb^{(0)\top}\Fb^{(0)}\right)\Db^{1/2}\Qb^{(0)\top}\hat{\Qb}^{(0)}+\hat{\Qb}^{(0)\top}\Qb^{(0)}\Db^{1/2}\left(\frac{1}{T_0}\Fb^{(0)\top}\Eb^{(0)}\right)\hat{\Qb}^{(0)}\\
&+\hat{\Qb}^{(0)\top}\left(\frac{1}{T_0}\Eb^{(0)\top}\Fb^{(0)}\right)\Db^{1/2}\Qb^{(0)\top}\hat{\Qb}^{(0)}+\hat{\Qb}^{(0)\top}\left(\frac{1}{T_0}\Eb^{(0)\top}\Eb^{(0)}\right)\hat{\Qb}^{(0)}.
 \end{aligned}
\end{equation*}
 By the weak dependence assumptions on $\bm{f}_t$ and $\bepsilon_t$ in Assumption \ref{corr}, the last three terms are asymptotically negligible compared to the first term. Combining with Assumption \ref{factors} which implies $\Fb^{(0)\top}\Fb^{(0)}/{T_0} = \bSigma_F^{(0)}+o_p(1)$, 
 we can obtain 
 
 $$\hat{\Mb}  \approx \hat{\Qb}^{(0)\top}\Qb^{(0)}\Db^{1/2}\left(\frac{1}{T_0}\Fb^{(0)\top}\Fb^{(0)}\right)\Db^{1/2}\Qb^{(0)\top}\hat{\Qb}^{(0)}\approx \Db^{1/2}\bSigma_F\Db^{1/2}.$$ 
Therefore, $\hat{\Mb}_{ii}\approx \left(\bSigma_{F}\right)_{ii} \Db_{ii}$, where $\left(\bSigma_{F}\right)_{ii}$ represent the $i$-th diagonal element of $\bSigma_F$. 
Recall that in Assumption \ref{loadings}.(2), the true factor strength $\alpha_i$ is defined as 
$$\left(\Db_w\right)_{ii}=C_i N^{\alpha_i},\quad i \in [s],$$
where $C_i$ is a constant that varies with $i$. Then we have
$$\tilde{\alpha}_i=\frac{\log \left(\Db_w\right)_{ii}}{\log(N)}=\alpha_i+\frac{\log(C_i)}{\log N}, \quad i\in[s].$$
As the dimension $N$  tends to infinity, we have ${\log(C_i)}/{\log N} \rightarrow 0$ and then $\tilde{\alpha}_i\rightarrow \alpha_i$. In particular, if $C_i=1$, we have $\tilde{\alpha}_i=\alpha_i$. By using the diagonal entries of $\hat{\Mb}$,  we can directly derive the estimator for $\Db$, and the estimators of  the factor strengths are denoted as $\hat{\alpha}_i^\ast,\ i=1,\ldots,r_0$, where
\begin{equation}\label{fs}
\hat{\alpha}_i^\ast=\frac{\log\left(\hat{\Mb}_{ii}\right)}{\log N}.\end{equation}
For the estimators
$\hat{\alpha}_i^\ast,\ i=1,\ldots,r_0$, the first  
$(r_0-s)$ ones correspond to strong factor strengths and the last $s$ ones correspond to weak factor strengths. The following theorem establishes the consistency for the estimators of the factor strengths $\hat{\alpha}_i^\ast$ in (\ref{fs}) under mild conditions.

\begin{theorem}\label{th5}
Under Assumptions \ref{factors} to \ref{corr}, we have
$$\hat{\alpha}_i^\ast-\alpha_i^*=O_p\left(\frac{1}{\log N}\right),\ i=1,\ldots,r_0,$$
where $\alpha_i^\ast=1$ for $i\in[r_0-s]$ and $\alpha^\ast_{r_0-s+i}=\alpha_i$ for $i\in[s]$.
\end{theorem}

In the following we introduce another method to estimate the factor strenghths, which is associated with the TransPCA procedure with  auxiliary datasets.
In the Step 2 of Algorithm \ref{algorithm:1}, we obtain the estimator $\hat{\Qb}^{(0)}_w$ of $\Qb_w^{(0)}$. Due to  the fact that $\hat{\Qb}_w^{(0)\top}\Qb_s^{(0)} \approx \bm{0}$,  we can also get the following approximation
\begin{equation*}
\begin{aligned}
\hat{\Mb}_w &=\hat{\Qb}_w^{(0)\top}\hat{\bSigma}^{(0)}\hat{\Qb}_w^{(0)}\\
&\approx \left(\begin{array}{cc}\hat{\Qb}_w^{(0)\top} \Qb_w^{(0)}\Db_w^{1/2} &\bm{0}\end{array}\right) \left(\frac{1}{T_0} \Fb^{(0)\top}\Fb^{(0)}\right) \left(
 \begin{array}{c}
\Db_w^{1/2}\Qb_w^{(0)\top}\hat{\Qb}_w^{(0)}  \\
   \bm{0}
  \end{array}
  \right) \\
  &\approx \Db_w^{1/2}\bSigma_F^{\ast}\Db_w^{1/2},
\end{aligned}
\end{equation*}
where $\Fb_w^{(0)\top}\Fb_w^{(0)}/T_0\overset{P}{\rightarrow}\bSigma_F^\ast$.
Similarly, the factor strength $\alpha_{i}$ can be estimated by
\begin{equation}\label{fs-trans}
\hat{\alpha}_{i} = \frac{\log\left((\hat{\Mb}_{w})_{ii}\right)}{\log N}, \ i \in [s],
\end{equation}
and then let $\hat{\Db}_w=\text{diag}\left(N^{\hat{\alpha}_1},\cdots,N^{\hat{\alpha}_s}\right)$. The estimator of the  loading matrix corresponding to weak factors can naturally be constructed as 
$\hat{\bLambda}^{(0)}_w = \hat{\Qb}^{(0)}_w\hat{\Db}_w^{1/2}$.
The following
Theorem \ref{th6} shows the consistency of the estimators of  factor strengths defined in (\ref{fs-trans}).

\begin{theorem}\label{th6}
Under Assumptions \ref{factors} to  \ref{NT-relation}, we further assume that  $\varepsilon=o\left({1}/{\sqrt{N^{1-\alpha_s}}}\right)$. As $N$, $ \underset{k\in\{0\}\cup \mathcal{A}}{\min}T_k \rightarrow \infty$, we have
$$\hat{\alpha}_i-\alpha_i=O_p\left(\frac{1}{\log N}\right),\ i\in [s].$$
\end{theorem}
Theorem \ref{th5} and \ref{th6}  assert that both the estimators $\hat{\alpha}_{i}$ and $\hat{\alpha}_i^\ast$ of the factor strength  converge to the true factor strength with a rate of $1/\log N$. By Theorem \ref{th5} and \ref{th6}, we have the option to estimate the factor strengths either solely based on the target dataset or  leveraging the auxiliary dataset. It is reasonable for both methods to yield consistent estimators with the same convergence rate, as mentioned earlier, our TransPCA method essentially improves the estimation accuracy by increasing the effective sample size. Theorem \ref{th5} indicates that the estimation accuracy  of the factor strengths does not depend on the sample size of the target dataset $T_0$, and therefore the TransPCA method  would not be able to improve  estimation accuracy of the  factor strengths.

\section{Useful Dataset Selection}\label{sec:data selection}
In Section \ref{sec:Methodology}, we assume that all informative datasets $\big\{\Xb^{(k)}, k \in \mathcal{A} \big\}$ are known in advance and the informative level is measured by equation (\ref{similarity}). Corollary \ref{cor:1}, Theorem \ref{th2} and Theorem \ref{th6}  all assume that $\varepsilon=o\left({1}/{\sqrt{N^{1-\alpha_s}}}\right)$, which ensures that the oracle TransPCA estimators are consistent.
However, in practice, we do not know the index set $\mathcal{A}$ in advance. Therefore, we need to select informative datasets from a huge amount of source panels to avoid negative transfer. 

In the first step of the oracle TransPCA procedure,
$\hat{\Qb}_w^{(0)}$ is composed of the leading $s$ eigenvectors of $\hat{\Pb}^w$, which can also be viewed as the solution of the 
following optimization problem
$$\hat{\Qb}_w^{(0)}=\underset{\Qb^\top\Qb=\Ib_s, \ \Qb \in\RR^{N\times s}}{\arg\max}\frac{1}{T}\sum_{k\in\{0\}\cup\mathcal{A}}T_k\text{tr}\left(\Pb_{\hat{\Qb}^{(k)}}\Pb_{\Qb}\right).$$
Note that under Assumption \ref{Identifiability},
\begin{equation}\label{equ:informative data}
s-\frac{\varepsilon^2}{2}\leq \text{tr}\left(\Pb_{\Qb^{(k)}}\Pb_{\Qb_w^{(0)}}\right)\leq s, \ k\in \mathcal{A},
\end{equation}
the informative level of the auxiliary datasets can be measured by $\text{tr}\left(\Pb_{\Qb^{(k)}}\Pb_{\Qb_w^{(0)}}\right)$, which motivates us to consider the following rectified problem
\begin{equation}\label{equ:selec data}
\hat{\Qb}_{w}^{(0)}(\tau)=\underset{\Qb^\top\Qb=\Ib_s, \ \Qb \in\RR^{N\times s}}{\arg\max}\frac{1}{\tilde{T}}
\left(T_0\text{tr}\left(\Pb_{\hat{\Qb}^{(0)}}\Pb_\Qb\right)+\sum_{k\in[K]}T_k\max\Big\{\text{tr}\left(\Pb_{\hat{\Qb}^{(k)}}\Pb_{\Qb}\right),\tau\Big\}\right),
\end{equation}
where $\tilde{T}=\sum_{k\in \{0\}\cup [K]} T_k$, $\tau \in [0,s]$ is a threshold parameter. If $\tau=0$, it indicates that all auxiliary datasets are integrated together to estimate $\Qb_w^{(0)}$  without any selection of datasets. If $\tau=s$, it is equivalent to performing PCA solely based on the target dataset.

Note that the optimization problem in \eqref{equ:selec data} involves combinatorial non-convex optimization, which makes it  difficult and computationally expensive to obtain the global maximizer. As an alternative, we adopt the cyclic coordinate descent algorithm to numerically search for the local maximum of (\ref{equ:selec data}). Specifically in the $m$-th iteration with a given estimator $\Pb_{\hat{\Qb}^{(0)}_{w,m-1}}$, we select the informative datasets according to 
\begin{equation}\label{equ:selectA}
\hat{\mathcal{A}}^{(m)}=\bigg\{k\in[K]: \text{tr}\left(\Pb_{\hat{\Qb}^{(k)}}\Pb_{\hat{\Qb}^{(0)}_{w,m}}\right)\geq \tau\bigg\}.
\end{equation}
Similarly,  given $\hat{\mathcal{A}}^{(m)}$, we can update $\hat{\Qb}^{(0)}_{w,m}$ as the leading $s$ eigenvectors of $\sum_{k\in \hat{\mathcal{A}}^{(m)}}T_k\Pb_{\hat{\Qb}^{(k)}}/\left(\sum_{k\in \hat{\mathcal{A}}^{(m)}}T_k\right)$.
We summarize the details in Algorithm \ref{alg2} for better illustration.
\begin{algorithm}
\caption{Non-oracle TransPCA Procedure}\label{alg2}
\begin{algorithmic}[1]
\REQUIRE Datasets $\big\{\Xb^{(k)}, k\in \{0\}\cup [K]\big\}$; $\Db_w$; initial estimator $\hat{\Qb}_{w,0}^{(0)}$; threshold parameter $\tau$; maximum iteration $m_{\max}$.
\WHILE{$m\le m_{\max}$ or $\hat{\mathcal{A}}^{(m-1)}
 \neq \hat{\mathcal{A}}^{(m)}$}
    \STATE $m\xleftarrow{} m+1$.
    \STATE Select the informative datasets $\hat{\mathcal{A}}^{(m)}$ based on \eqref{equ:selectA}.
    \STATE Update $\hat{\Qb}^{(0)}_{w,m}$ as the leading $s$ eigenvectors of $\sum_{k\in \hat{\mathcal{A}}^{(m)}}T_k\Pb_{\hat{\Qb}^{(k)}}/\left(\sum_{k\in \hat{\mathcal{A}}^{(m)}}T_k\right)$.
    \ENDWHILE
\ENSURE Selected informative datasets $\hat{\mathcal{A}}$ and the corresponding estimator of weak loading matrix $\hat{\Qb}^{(0)}_{w}$. 
\end{algorithmic}
\end{algorithm}


Theoretically,
to guarantee the consistency of the dataset selection procedure, we need an additional assumption to ensure a sufficient distinction between transferable and non-transferable datasets.

\begin{assumption}\label{asmp:non-informative}
For any $k\in [K]\setminus \mathcal{A}$, we assume that there exists $h_\tau>0$ such that $\text{tr}\left(\Pb_{\Qb^{(k)}}\Pb_{\Qb_w^{(0)}}\right)\leq s-h_\tau$.
\end{assumption}

Assumption \ref{asmp:non-informative} ensures the separability between informative and non-informative datasets, ensuring informative source models  more similar to the target model.

\begin{theorem}\label{th7}
Under Assumptions \ref{factors} to \ref{NT-relation}.(1) and Assumption \ref{asmp:non-informative}, if $\sqrt{N^{1-\alpha_1}}\delta+\underset{k\in[K]\setminus \mathcal{A}}{\max}\ T_k^{-1/2}=o\left(h_\tau\right)$ as $N$, $ \underset{k\in\{0\}\cup [K]}{\min}T_k \rightarrow \infty$ and $\varepsilon \rightarrow 0$, for any $\tau\in[s-c_1 h_\tau, s-c_2 h_\tau]$, where $c_1$ and $c_2$ are constants satisfying $0<c_2<c_1<1$,  the weak loading matrix estimator $\hat{\Qb}_w^{(0)}$ is a local maximum of (\ref{equ:selec data}) with probability tending to 1.
\end{theorem}
In Theorem \ref{th7}, we do not require $h_\tau$ to be a constant; instead, we allow it to tend to zero, provided that $\sqrt{N^{1-\alpha_1}}\delta+\underset{k\in[K]\setminus \mathcal{A}}{\max}\ T_k^{-1/2}=o\left(h_\tau\right)$. 
As $N$, $ \underset{k\in\{0\}\cup [K]}{\min}T_k \rightarrow \infty$, $\varepsilon\rightarrow 0$, it always holds that $\text{tr}\left(\Pb_{\Qb^{(k)}}\Pb_{\Qb_w^{(0)}}\right)< \text{tr}\left(\Pb_{\Qb^{(l)}}\Pb_{\Qb_w^{(0)}}\right)$, $k\in [K]\setminus \mathcal{A}$, $l\in \mathcal{A}$, enabling us to select datasets by choosing an appropriate parameter $\tau$ that meets the conditions of Theorem \ref{th7}. In our simulation studies and real applications, the optimal parameter $\tau$ can be selected using $k$-fold cross-validation, based on the criterion of minimizing the Mean Squared Error (MSE) of the target dataset. Moreover, we also need to find a good initial estimator and utilize cyclical gradient descent method to find the local maximum. We can use $\hat{\Qb}^{(0)}$, obtained through performing PCA on the target dataset, as the initial estimator $\hat{\Qb}_{w,0}^{(0)}$ for Algorithm \ref{alg2}.

\section{Simulation Study}\label{sec:simulation}
In this section, we conduct simulation studies to evaluate the finite sample performance of our proposed methods. In Section \ref{DGP}, we introduce the data generating mechanism for the target and auxiliary datasets.  We evaluate the finite sample performance of the Oracle TransPCA procedure in terms of estimating the loading spaces and common components of the target model in Section \ref{sec:simloading}. In Section \ref{sec:fs}, we  investigate the finite sample performances of the proposed two estimators of factor strengths. At last, in Section  \ref{sec:fn}, we evaluate the finite sample performances of the proposed ED estimator of the number of weak factors.

\subsection{Data generation}\label{DGP}

 We generate the simulated target and source datasets as follows:
\begin{itemize}
    \item Factor loadings $\bLambda^{(k)}$, $k\in \{0\} \cup [K]$, \vspace{0.5em}
    \\ 
    ($i$) for the target dataset, generate an $N\times s$ column orthogonal matrix $\Qb_w^{(0)}$ and an $N\times (r_0-s)$  column orthogonal matrix $\Qb_s^{(0)}$ from the orthogonal complement of  $\mathcal{M}\left(\Qb_w^{(0)}\right)^\perp$. Then we let $\bLambda^{(0)}=\left(\begin{array}{cc}\Qb_w^{(0)}\Db_w^{1/2}& \sqrt{N}\Qb_s^{(0)}\end{array}\right)$; \vspace{0.5em}
    \\
    ($ii$) for the informative auxiliary dataset $k \in\mathcal{A}$, we first acquire an orthogonal matrix $\Ub^{(k)}$ close to $\Ib_N$ by taking QR decomposition  of the perturbed identity matrix $ \Ib_N+\mathcal{E}^{(k)}=\Ub^{(k)}\Rb^{(k)}$, where the elements of $\mathcal{E}^{(k)}$ are drawn
independently from $\mathcal{N} (0, {\varepsilon^2}/{N})$, then we let  $\Qb_\mathcal{T}^{(k)} = \Ub^{(k)} \Qb_w^{(0)}$. Generate  an $N \times (r_k - s)$-dimensional column orthogonal matrix $\Qb_{\mathcal{T}^c}^{(k)}$ from $\mathcal{M}(\Qb_\mathcal{T}^{(k)})^\perp$. Then we obtain the loading matrix $\bLambda^{(k)} = \sqrt{N}\left(\begin{array}{cc}\Qb_\mathcal{T}^{(k)}& \Qb_{\mathcal{T}^c}^{(k)}\end{array}\right)$; \vspace{0.5em}
\\
($iii$) for the non-informative dataset $k\in [K]\setminus \mathcal{A}$, we first generate an $N\times s$ column orthogonal matrix $\tilde{\Qb}_w^{(0)}$ whose spanned space is different from that spanned by the columns of $\Qb_w^{(0)}$. Then we let $\Qb_\mathcal{T}^{(k)} = \Ub^{(k)} \tilde{\Qb}_w^{(0)}$ and $\bLambda^{(k)} = \sqrt{N}\left(\begin{array}{cc}\Qb_\mathcal{T}^{(k)}& \Qb_{\mathcal{T}^c}^{(k)}\end{array}\right)$, where $\Ub^{(k)}$ and $\Qb_{\mathcal{T}^c}^{(k)}$  are  generated similar to the procedure in ($ii$).

    \item Factors and idiosyncratic errors, for each $ k\in \{0\} \cup [K]$,  let
    $$
\bm{f}_t^{(k)} = \varphi \bm{f}_{t-1}^{(k)} + \sqrt{1-\varphi^2}\bm{u}_t^{(k)},  \ \ \
\bepsilon_t^{(k)} = \phi \bepsilon_{t-1}^{(k)} + \sqrt{1-\phi^2}\bm{v}_t^{(k)}, 
$$ 
where $\varphi$ and $\phi$ control the temporal correlations, $\bu_t^{(k)}$ and 
$\bv_t^{(k)}$
are $i.i.d.$  random samples from multivariate Gaussian distributions $\mathcal{N}(\bm{0},\Ib_{r_k}$) and $\mathcal{N}(\bm{0},\Ib_{N}$), respectively.
 \end{itemize}

The following scenarios are considered: ($\bm{a}$) all $K$ auxiliary datasets are informative; ($\bm{b}$) $K/2$ auxiliary datasets are informative. We set $r_0 = 3$, $\varepsilon = 0.1$, $s = 2$, $r_k = 4$ for all $k \in [K]$,
and the weak factor strengths $(\alpha_{1},\alpha_{2})=(0.7,0.6)$, i.e., $\Db_w=\text{diag}\left(N^{0.7}, N^{0.6}\right)$.  We consider all combinations of dimension $N \in \{50, 100\}$,  $T_0 \in \{50, 100\}$,  $T_k \in \{200,300,400\}$ for $ k \in [K]$ and the number of auxiliary datasets $K \in \{4, 8\}$.


\begin{table}[!h]
    \centering
    	\scalebox{1}{ \renewcommand{\arraystretch}{1}
    \begin{tabular}{ccccccccccc} 
    \toprule
       $K$ & $N$ & $T_0$ & $T_k$ & $\mathcal{D}$(Oracle) & $\mathcal{D}$(Target-only) & MSE(Oracle) & MSE(Target-only)\\
        \hline

       \multirow{12}{*}{4} & \multirow{6}{*}{50} & \multirow{3}{*}{50} & 200 & 0.087(0.011) & 0.249(0.024) & 0.082(0.009) & 0.128(0.011)\\
        &  &  & 300 & 0.085(0.011) & 0.249(0.025) & 0.081(0.009) & 0.128(0.011) \\
        &  &  & 400 & 0.083(0.011) & 0.250(0.024) & 0.080(0.008) & 0.128(0.011) \\

        &  & \multirow{3}{*}{100} & 200 & 0.066(0.007) & 0.175(0.015) & 0.071(0.006) & 0.094(0.007) \\
        &  &  & 300 & 0.063(0.007) & 0.177(0.015) & 0.071(0.005) & 0.094(0.006) \\
        &  &  & 400 & 0.061(0.007) & 0.176(0.015) & 0.071(0.005) & 0.094(0.007) \\
        & \multirow{6}{*}{100} & \multirow{3}{*}{50} & 200 & 0.088(0.009) & 0.277(0.022) & 0.051(0.004) & 0.096(0.006) \\
        &  &  & 300 & 0.086(0.009) & 0.277(0.022) & 0.051(0.004) & 0.097(0.006) \\
        &  &  & 400 & 0.085(0.010) & 0.278(0.023) & 0.051(0.005) & 0.096(0.007) \\
        &  & \multirow{3}{*}{100} & 200 & 0.067(0.005) & 0.197(0.014) & 0.041(0.003) & 0.063(0.004) \\
        &  &  & 300 & 0.064(0.005) & 0.197(0.014) & 0.041(0.003) & 0.063(0.004) \\
        &  &  & 400 & 0.062(0.006) & 0.196(0.013) & 0.041(0.003) & 0.063(0.004) \\

\hline
\multirow{12}{*}{8} & \multirow{6}{*}{50} & \multirow{3}{*}{50} & 200 & 0.083(0.011) & 0.249(0.025) & 0.080(0.008) & 0.128(0.012)\\
        &  &  & 300 & 0.082(0.010) & 0.248(0.025) & 0.081(0.008) & 0.129(0.011) \\
        &  &  & 400 & 0.083(0.012) & 0.250(0.026) & 0.080(0.008) & 0.128(0.011) \\
        &  & \multirow{3}{*}{100} & 200 & 0.061(0.007) & 0.175(0.014) & 0.071(0.005) & 0.094(0.007) \\
        &  &  & 300 & 0.060(0.007) & 0.176(0.015) & 0.071(0.005) & 0.094(0.007) \\
        &  &  & 400 & 0.059(0.007) & 0.176(0.014) & 0.070(0.006) & 0.094(0.007) \\
        & \multirow{6}{*}{100} & \multirow{3}{*}{50} & 200 & 0.085(0.010) & 0.277(0.024) & 0.051(0.005) & 0.096(0.006) \\
        &  &  & 300 & 0.084(0.010) & 0.275(0.023) & 0.051(0.004) & 0.096(0.006) \\
        &  &  & 400 & 0.083(0.010) & 0.275(0.023) & 0.051(0.005) & 0.096(0.007) \\
        &  & \multirow{3}{*}{100} & 200 & 0.062(0.005) & 0.196(0.014) & 0.040(0.003) & 0.063(0.004) \\
        &  &  & 300 & 0.060(0.005) & 0.197(0.014) & 0.040(0.003) & 0.063(0.004) \\
        &  &  & 400 & 0.059(0.006) & 0.198(0.013) & 0.040(0.025) & 0.063(0.090) \\
      \hline  

        \toprule
\end{tabular}}
    \caption{ Averaged estimation errors and standard errors (in parentheses) of $\mathcal{D}\left(\mathcal{M}\left(\hat{\bLambda}^{(0)}\right),\mathcal{M}\left(\bLambda^{(0)}\right)\right)$ for Scenario ($\bm{a}$), where all $K$ auxiliary datasets are informative,
under different settings over 500 replications. ``Oracle" stands for the Oracle TransPCA method;  ``Target-only" stands for the PCA method performed solely on target dataset.}
    \label{tab:1}
\end{table}

\subsection{ Estimation of loading spaces and common components}\label{sec:simloading}
In this section, we compare the performance of the oracle TransPCA method (where the informative datasets are known and only the informative datasets are used), non-oracle TransPCA method (where the informative datasets are unknown and selected by Algorithm \ref{alg2}), target-only PCA method (where only the target dataset is used for PCA),  and blind TransPCA method (where all source datasets are blindly combined to perform the TransPCA procedure), in terms of estimating loading spaces and  common components.
\begin{figure}[!h]
\begin{minipage}[t]{0.48\linewidth}	
\centering
\includegraphics[width=1\textwidth]{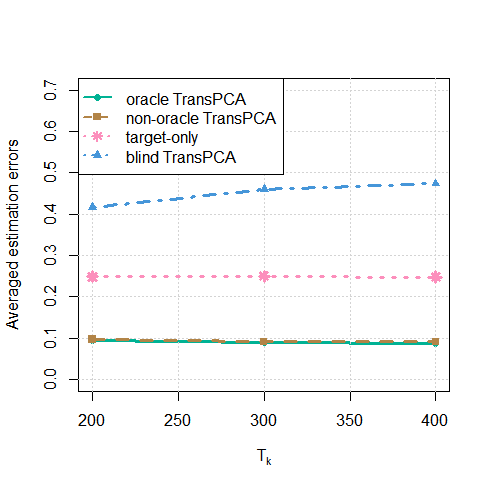}	
\end{minipage}
\hfill
\begin{minipage}[t]{0.48\linewidth}
\centering
\includegraphics[width=1\textwidth]{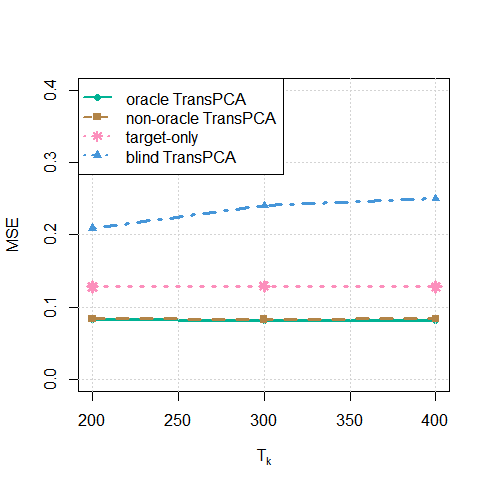}
\end{minipage}
\caption{Averaged estimation errors of $\mathcal{D}\left(\mathcal{M}\left(\hat{\bLambda}^{(0)}\right),\mathcal{M}\left(\bLambda^{(0)}\right)\right)$ (left figure) and MSE (right figure) for Scenario ($\bm{b}$) using different methods with $T_k=(200, 300, 400)$ when $K=4$, $N=50$, $T_0=50$.}\label{fig:1}
\end{figure}

To assess the empirical performance of the proposed estimators for loading spaces, we adopt a metric between linear spaces which was also utilized in \cite{he2022large}. For two column-orthogonal matrices $\Ob_1$ and $\Ob_2$ of sizes $N\times p_1$ and $N\times p_2$ ($\max\{p_1,p_2\}\leq N$), define 
$$\mathcal{D}\left(\mathcal{M}(\Ob_1),\mathcal{M}(\Ob_2)\right)=\left(1-\frac{\text{tr}(\Ob_1\Ob_1^\top\Ob_2\Ob_2^\top)}{\min\{p_1,p_2\}}\right)^{1/2}.$$
The distance measure $\mathcal{D}\left(\mathcal{M}(\Ob_1),\mathcal{M}(\Ob_2)\right)$ is a quantity between 0 and 1. It is equal to 0 if $\mathcal{M}(\Ob_1) \subset \mathcal{M}(\Ob_2) $ or $\mathcal{M}(\Ob_2) \subset \mathcal{M}(\Ob_1) $, and is equal to 1 if $\mathcal{M}(\Ob_1) $ and $\mathcal{M}(\Ob_2) $ are orthogonal. We measure the estimation error of the common components by using the Mean Squared Error, i.e.,
$$
\text{MSE} = \frac{1}{NT_0 } \left\|\hat{\Fb}^{(0)}\hat{\bLambda}^{(0)\top}-\Fb^{(0)}\bLambda^{(0)\top}\right\|_F^2.
$$

Firstly, we assume that all $K$ datasets are  informative. Table \ref{tab:1} shows the averaged estimation errors of $\mathcal{D}\left(\mathcal{M}\left(\hat{\bLambda}^{(0)}\right),\mathcal{M}\left(\bLambda^{(0)}\right)\right)$ and the MSE values of the common components over 500 replications. From Table \ref{tab:1}, it can be seen that our TransPCA method outperforms the traditional PCA performed solely to the target dataset,  with both the loading matrix estimation errors and the MSEs being smaller across all settings. 
When $N$ and $T_0$ are fixed, as the sample size of the auxiliary datasets increases, our TransPCA method achieves higher accuracy in terms of estimating loading matrices. This further validates that our method inherently improves estimation accuracy by increasing the total effective sample size.

Next, we consider the scenario where some of the panel datasets are non-informative, necessitating the useful dataset selection step to avoid negative transfer.  Figure \ref{fig:1}  shows the averaged estimation errors and MSEs for  different methods  as $T_k$ changes, with $K=4$, $N=50$ and $T_0=50$. Tables \ref{tab:2}  and \ref{tab:4} show the averaged estimation errors of $\mathcal{D}\left(\mathcal{M}\left(\hat{\bLambda}^{(0)}\right),\mathcal{M}\left(\bLambda^{(0)}\right)\right)$ and MSEs with standard errors in parentheses under Scenario ($\bm{b}$) for different settings. 
The same conclusions can be drawn from Tables \ref{tab:2} and \ref{tab:4}  as those drawn from Table \ref{tab:1} that oracle TransPCA method performs better than target-only method. The non-oracle TransPCA method performs comparably with oracle TransPCA, which  indicates that the dataset selection procedure proposed in Section \ref{sec:data selection} effectively identifies informative datasets as if they were known a priori. In addition, we can see that blind TransPCA performs the worst among all methods, suggesting that when non-informative datasets are present, blindly combining all datasets  may lead to negative transfer, deteriorating the estimation accuracy for the target model.

We further evaluate the accuracy of the proposed dataset selection method using the following metrics: True Positive Rate (TPR), True Negative Rate (TNR) and Precision. We label the $K$ auxiliary datasets, assigning the class label ``positive'' to informative datasets and ``negative'' to non-informative datasets.
{The results in Table \ref{tab:5} demonstrate that our dataset selection method almost accurately identifies informative datasets, with TPR, TNR and precision uniformly approaching to  1 across various settings.}

\begin{table}[!h]
    \centering
    \scalebox{1}{\renewcommand{\arraystretch}{1.1}
    \begin{tabular}{ccccccccccc}
    \toprule
       $K$ & $N$ & $T_0$ & $T_k$ & \text{Oracle } & Non-oracle & Target-only& Blind \\
        \hline
       \multirow{12}{*}{4} & \multirow{6}{*}{50} & \multirow{3}{*}{50} & 200 & 0.095(0.014)& 0.097(0.028)&0.249(0.023)&0.417(0.065)\\
        &  &  & 300&0.089(0.012)&0.091(0.018)&0.250(0.024)&0.460(0.073)\\
        &  &  & 400 &0.087(0.011)&0.091(0.043)&0.248(0.025)&0.474( 0.076) \\

        &  & \multirow{3}{*}{100} & 200 &0.076(0.006)&0.077(0.021)&0.176(0.015)&0.322(0.051)\\
        &  &  & 300 &0.069(0.006)&0.071(0.027)&0.176(0.016)&0.375(0.059)\\
        &  &  & 400 &0.066(0.007)& 0.068(0.029)&0.177(0.014)&0.407(0.065)\\
        & \multirow{6}{*}{100} & \multirow{3}{*}{50} & 200 & 0.095(0.009)&0.095(0.009)&0.278(0.023)&0.382(0.052)\\
        &  &  & 300 &0.090(0.009)&0.091(0.024)&0.276(0.024)&0.429(0.061)\\
        &  &  & 400 &0.088(0.009)&0.088(0.009)&0.277(0.023)&0.451(0.063) \\
        &  & \multirow{3}{*}{100} & 200 &0.077(0.005)&0.077(0.005)&0.196(0.013)&0.284(0.045) \\
        &  &  & 300 &0.071(0.005)&0.071(0.005)&0.197(0.014)& 0.337(0.046)\\
        &  &  & 400 & 0.067(0.005)& 0.067(0.005)& 0.197(0.013)&0.370(0.057)\\

\hline
\multirow{12}{*}{8} & \multirow{6}{*}{50} & \multirow{3}{*}{50} & 200 &0.086(0.011)&0.086(0.012)&0.248(0.023)&0.443(0.055)\\
        &  &  & 300&0.085(0.011)&0.086(0.012)&0.250(0.025)&0.474(0.057)\\
        &  &  & 400 &0.084(0.011)&0.084(0.012)&0.250(0.025)&0.484(0.065) \\
        &  & \multirow{3}{*}{100} & 200 &0.066(0.006)&0.066(0.006)&0.177(0.015)&0.372(0.047) \\
        &  &  & 300 &0.063(0.006)&0.063(0.007)&0.176(0.015)&0.409(0.053) \\
        &  &  & 400&0.062(0.007)& 0.062(0.008)&0.176(0.015)&0.433(0.056)\\
        & \multirow{6}{*}{100} & \multirow{3}{*}{50} & 200& 0.088(0.010)&0.088(0.010)&0.275(0.023)&0.429(0.050) \\
        &  &  & 300&0.086(0.010)&0.086(0.010)&0.276(0.022)&0.458(0.052)\\
        &  &  & 400&0.085(0.010)&0.085(0.010)&0.275(0.022)&0.474(0.052)\\
        &  & \multirow{3}{*}{100} & 200 &0.067(0.005)&0.067(0.005)&0.196(0.013)&0.352(0.047) \\
        &  &  & 300&0.063(0.006)&0.063(0.006)&0.197(0.014)&0.395(0.046)\\
        &  &  & 400 &0.062(0.005)&0.062(0.005)&0.196(0.013)&0.418(0.050)\\
      \hline  
         \toprule
    \end{tabular}}
    \caption{
    Averaged estimation errors and standard errors (in parentheses) of $\mathcal{D}\left(\mathcal{M}\left(\hat{\bLambda}^{(0)}\right),\mathcal{M}\left(\bLambda^{(0)}\right)\right)$ for Scenario ($\bm{b}$), where $K/2$ auxiliary datasets are informative,
under different settings over 500 replications. ``Oracle" stands for the Oracle TransPCA method; ``Non-Oracle" stands for the  TransPCA method with useful dataset selection capability; ``Target-only" stands for the PCA method performed solely on target dataset; ``Blind" stands for the TransPCA method with all auxiliary datasets treated as useful.
  }
    \label{tab:2}
\end{table}

\begin{table}[!h]
    \centering
    \scalebox{1}{\renewcommand{\arraystretch}{1}
    \begin{tabular}{ccccccccccc}
    \toprule
       $K$ & $N$ & $T_0$ & $T_k$ & Oracle & Non-oracle & Target-only & Blind \\
        \hline
       \multirow{12}{*}{4} & \multirow{6}{*}{50} & \multirow{3}{*}{50} & 200 &0.083(0.008)&0.084(0.011)&0.128(0.010)&0.209(0.049)\\
        &  &  & 300 &0.082(0.009)&0.083(0.010)&0.129(0.011)&0.240(0.063)\\
        &  &  & 400&0.081(0.009)&0.084(0.024)&0.128(0.011)&0.250(0.062)&\\

        &  & \multirow{3}{*}{100} & 200 &0.073(0.006)&0.073(0.009)&0.094(0.006)&0.150(0.031) \\
        &  &  & 300&0.072(0.006)&0.072(0.014)&0.094(0.007)& 0.178(0.040)\\
        &  &  & 400 &0.071(0.005)& 0.072(0.015)&0.094(0.006)&0.198(0.047)\\
        & \multirow{6}{*}{100} & \multirow{3}{*}{50} & 200 & 0.052(0.005)&0.052(0.005)&0.097(0.007)&0.134(0.031)\\
        &  &  & 300 & 0.052(0.005)&0.052(0.010)&0.096(0.007)&0.160(0.043)\\
        &  &  & 400 &0.051(0.005)&0.051(0.005)&0.097(0.007)&0.172(0.041)\\
        &  & \multirow{3}{*}{100} & 200 &0.042(0.003)&0.042(0.003)&0.063(0.004)&0.089(0.017) \\
        &  &  & 300 &0.042(0.003)&0.042(0.003)&0.063(0.004)&0.109(0.022) \\
        &  &  & 400&0.041(0.003)&0.041(0.003)&0.063(0.004)&0.123(0.030) \\

\hline
\multirow{12}{*}{8} & \multirow{6}{*}{50} & \multirow{3}{*}{50} & 200&0.082(0.008)&0.082(0.008)&0.128(0.011)&0.227(0.044)\\
        &  &  & 300&0.081(0.008)&0.081(0.008)& 0.128(0.011)&0.247(0.049)\\
        &  &  & 400&0.081(0.009)&0.081(0.009)&0.128(0.011)& 0.255(0.055)\\
        &  & \multirow{3}{*}{100} & 200&0.071(0.005)&0.071(0.005)&0.094(0.007)&0.175(0.031)\\
        &  &  & 300 &0.071(0.005)&0.071(0.005)&0.094(0.007)&0.198(0.036)\\
        &  &  & 400 &0.071(0.005)&0.071(0.005)&0.094(0.007)&0.214(0.041) \\
        & \multirow{6}{*}{100} & \multirow{3}{*}{50} & 200 &0.051(0.005)& 0.051(0.005)&0.096(0.008)&0.158(0.032)\\
        &  &  & 300 & 0.051(0.005)&0.051(0.005)&0.096(0.007)&0.173(0.034)\\
        &  &  & 400&0.051(0.005)&0.051(0.005)&0.096(0.007)&0.181(0.034) \\
        &  & \multirow{3}{*}{100} & 200 &0.041(0.003)&0.041(0.003)&0.063(0.004)&0.115(0.022) \\
        &  &  & 300&0.041(0.003)&0.041(0.003)&0.063(0.004)&0.134(0.024)\\
        &  &  & 400& 0.041(0.003)& 0.041(0.003)&0.063(0.004)&0.146(0.029)\\
      \hline  
         \toprule
    \end{tabular}}
    \caption{Mean squared error and its standard deviation under Scenario ($\bm{b}$) over 500 replications. ``Oracle" stands for the Oracle TransPCA method; ``Non-Oracle" stands for the  TransPCA method with useful dataset selection capability; ``Target-only" stands for the PCA method performed solely on target dataset; ``Blind" stands for the TransPCA method with all auxiliary datasets treated as useful.}
    \label{tab:4}
\end{table}

\begin{table}[!h]
    \centering
    \scalebox{1}{
    \begin{tabular}{ccccccccccc}
    \toprule
       $K$ & $N$ & $T_0$ & $T_k$ & TPR & TNR & Precision \\
        \hline

      \multirow{12}{*}{4} & \multirow{6}{*}{50} & \multirow{3}{*}{50} & 200 & 1.000(0.000) & 0.992(0.074) &  0.996(0.041) \\
        &  &  & 300 & 1.000(0.000) & 0.994(0.050) &0.996(0.033) \\
        &  &  & 400 & 1.000(0.000) & 0.992(0.059) &0.994(0.039) \\

        &  & \multirow{3}{*}{100} & 200 & 1.000(0.000) & 0.994(0.059)  &0.996(0.034)  \\
        &  &  & 300 & 1.000(0.000) & 0.992(0.059)  & 0.994(0.039) \\
        &  &  & 400 & 1.000(0.000) & 0.994(0.059) & 0.996(0.034) \\
        & \multirow{6}{*}{100} & \multirow{3}{*}{50} & 200 & 1.000(0.000) & 1.000(0.000)  & 1.000(0.000)\\
        &  &  & 300 & 1.000(0.000) & 0.998(0.022) & 0.998(0.015)\\
        &  &  & 400 & 1.000(0.000) & 1.000(0.000) & 1.000(0.000)\\
        &  & \multirow{3}{*}{100} & 200 & 1.000(0.000) & 1.000(0.000) & 1.000(0.000)\\
        &  &  & 300 & 1.000(0.000) & 1.000(0.000) & 1.000(0.000)\\
        &  &  & 400 & 1.000(0.000) & 1.000(0.000) &1.000(0.000) \\

\hline
\multirow{12}{*}{8} & \multirow{6}{*}{50} & \multirow{3}{*}{50} & 200 & 1.000(0.000) & 0.998(0.025) & 0.998(0.020)\\
        &  &  & 300 & 1.000(0.000) & 0.996(0.029) & 0.996(0.024)\\
        &  &  & 400 & 1.000(0.000) & 0.994(0.037) & 0.996(0.029) \\
       &  & \multirow{3}{*}{100} & 200 & 1.000(0.000) & 0.998(0.025) &0.998(0.017) \\
        &  &  & 300 & 1.000(0.000) & 0.998(0.022) & 0.998(0.018) \\
        &  &  & 400 & 1.000(0.000) & 0.998(0.019) &0.998(0.015) \\
        & \multirow{6}{*}{100} & \multirow{3}{*}{50} & 200 & 1.000(0.000) & 1.000(0.000) &  1.000(0.000)\\
        &  &  & 300 & 1.000(0.000) & 1.000(0.000) & 1.000(0.000) \\
       &  &  & 400 & 1.000(0.000) & 1.000(0.000) & 1.000(0.000)\\
        &  & \multirow{3}{*}{100} & 200 & 1.000(0.000) & 1.000(0.000) & 1.000(0.000)\\
        &  &  & 300 & 1.000(0.000) & 1.000(0.000) & 1.000(0.000)\\
        &  &  & 400 & 1.000(0.000) & 1.000(0.000) & 1.000(0.000) \\
      \hline  
         \toprule
    \end{tabular}}
    \caption{ Accuracy of the proposed data selection method (with standard deviation in parentheses) under different metrics for Scenario ($\bm{b}$) over 500 replications. }
    \label{tab:5}
\end{table}

\subsection{Estimation of the factor strengths}\label{sec:fs}
Recall that we estimate the factor strength $\alpha_i$ for $i\in [s]$ using the method described in Section \ref{sec:Estimate FS}. For brevity, in this section we denote the estimator $\hat{\alpha}_i^\ast$  as ``target-only" and the estimator $\hat{\alpha}_{i}$ as TransPCA, without causing any confusion. Tables \ref{label6} shows the mean and standard deviation of factor strength estimators  over 500 repetitions under each setting.
As shown in  Table \ref{label6}, for both  TransPCA method and target-only methods, factor strength estimators exhibit good performance in all settings and their performances are comparable, which is also consistent with the conclusion obtained in Theorem \ref{th5} and \ref{th6}. As the  dimension $N$ increases, both the  estimators of factor strength become more accurate.
\begin{table}[!h]
    \centering
    \scalebox{1}{\renewcommand{\arraystretch}{1}
    \begin{tabular}{ccccccccccc}
    \toprule
       $K$ & $N$ & $T_0$ & $T_k$ & $\hat{\alpha}_1$ & $\hat{\alpha}_2$ & $\hat{\alpha}_{1}^\ast$ & $\hat{\alpha}_{2}^\ast$\\
        \hline

       \multirow{12}{*}{4} & \multirow{6}{*}{50} & \multirow{3}{*}{50} & 200 & 0.71(0.052) & 0.62(0.053) & 0.73(0.045) & 0.62(0.045)\\
        &  &  & 300 & 0.71(0.050) & 0.62(0.060) & 0.73(0.044) & 0.62(0.048) \\
        &  &  & 400 & 0.72(0.051) & 0.62(0.055) & 0.73(0.044) & 0.62(0.044) \\

        &  & \multirow{3}{*}{100} & 200 & 0.71(0.040) & 0.63(0.041) & 0.72(0.034) & 0.62(0.033) \\
        &  &  & 300 & 0.71(0.038) & 0.63(0.043) & 0.72(0.033) & 0.62(0.036) \\
        &  &  & 400 & 0.71(0.038) & 0.63(0.040) & 0.73(0.033) & 0.62(0.033) \\
        & \multirow{6}{*}{100} & \multirow{3}{*}{50} & 200 & 0.71(0.043) & 0.61(0.044) & 0.72(0.038) & 0.62(0.037) \\
        &  &  & 300 & 0.71(0.042) & 0.60(0.044) & 0.72(0.038) & 0.62(0.037) \\
        &  &  & 400 & 0.71(0.040) & 0.60(0.045) & 0.72(0.036) & 0.62(0.037) \\
        &  & \multirow{3}{*}{100} & 200 & 0.71(0.030) & 0.61(0.033) & 0.72(0.027) & 0.61(0.029) \\
        &  &  & 300 & 0.71(0.031) & 0.61(0.035) & 0.72(0.028) & 0.62(0.030) \\
        &  &  & 400 & 0.71(0.033) & 0.61(0.033) & 0.72(0.029) & 0.62(0.028) \\

\hline
\multirow{12}{*}{8} & \multirow{6}{*}{50} & \multirow{3}{*}{50} & 200 & 0.71(0.049) & 0.62(0.055) & 0.73(0.042) & 0.62(0.044)\\
        &  &  & 300 & 0.71(0.052) & 0.62(0.056) & 0.73(0.043) & 0.62(0.045) \\
        &  &  & 400 & 0.71(0.053) & 0.62(0.059) & 0.73(0.047) & 0.62(0.046) \\
        &  & \multirow{3}{*}{100} & 200 & 0.71(0.041) & 0.63(0.043) & 0.73(0.034) & 0.62(0.034) \\
        &  &  & 300 & 0.71(0.039) & 0.63(0.042) & 0.73(0.033) & 0.62(0.032) \\
        &  &  & 400 & 0.71(0.040) & 0.63(0.043) & 0.72(0.033) & 0.62(0.033) \\
        & \multirow{6}{*}{100} & \multirow{3}{*}{50} & 200 & 0.71(0.041) & 0.61(0.050) & 0.72(0.036) & 0.62(0.040) \\
        &  &  & 300 & 0.71(0.039) & 0.61(0.046) & 0.72(0.036) & 0.62(0.037) \\
        &  &  & 400 & 0.71(0.044) & 0.60(0.046) & 0.73(0.039) & 0.62(0.038) \\
        &  & \multirow{3}{*}{100} & 200 & 0.70(0.033) & 0.62(0.036) & 0.71(0.029) & 0.62(0.030) \\
        &  &  & 300 & 0.71(0.031) & 0.61(0.035) & 0.72(0.028) & 0.61(0.030) \\
        &  &  & 400 & 0.70(0.032) & 0.61(0.034) & 0.71(0.027) & 0.61(0.029) \\
      \hline  
         \toprule
    \end{tabular}}
    \caption{The mean and standard deviation (in parentheses) for the estimated factor strengths $\hat{\alpha}_1$ ($\hat{\alpha}_1^\ast$) and $\hat{\alpha}_2$ ($\hat{\alpha}_2^\ast$), with true factor strengths $(\alpha_1,\alpha_2)=(0.7,0.6)$,
    under  different settings and Scenario ($\bm{a}$)
 over 500 replications.}
    \label{label6}
\end{table}

\subsection{Estimating the number of weak factors }\label{sec:fn}
In this section, we evaluate the empirical performance of our TransED method for estimating the number of weak factors. Table \ref{tab:weak number} presents the frequencies of exact estimation and underestimation of the number of weak factors over 500 replications. It can be seen from Table \ref{tab:weak number}  that  as the dimension  $N$ and sample size $T_0$ increase, the proportion of exact estimation by TransED method tends to converge to 1, for both the cases with the true number of weak factors being  2 or 3,  which is consistent with our theoretical analysis.

\begin{table}[!h]
    \centering
    \scalebox{1}{\renewcommand{\arraystretch}{1}
    \begin{tabular}{ccccccccccc}
    \toprule
     
        $s$&$K$ & $N$ & $T_0$ & $T_k=200$ & $T_k=300$ & $T_k=400$ \\
        \hline
        \multirow{8}{*}{2}&\multirow{4}{*}{4} & \multirow{2}{*}{50} & 50 & 0.992(0.000)&0.984(0.000)&0.992(0.000) \\
        & &  & 100 & 0.988(0.000)&0.990(0.000)&0.994(0.000) \\
        & & \multirow{2}{*}{100} & 50 & 1.000(0.000)&0.996(0.000)&0.996(0.000) \\
        & &  & 100 & 1.000(0.000)&1.000(0.000)&0.994(0.000) \\
        
        & \multirow{4}{*}{8} & \multirow{2}{*}{50} & 50 & 1.000(0.000)&1.000(0.000)&1.000(0.000) \\
        & &  & 100 &1.000(0.000)&1.000(0.000)&1.000(0.000) \\
        & & \multirow{2}{*}{100} & 50 & 1.000(0.000)&1.000(0.000)&1.000(0.000) \\
        & &  & 100 & 1.000(0.000)&1.000(0.000)&1.000(0.000) \\
        \hline
       \multirow{8}{*}{3}&\multirow{4}{*}{4} & \multirow{2}{*}{50} & 50 & 0.998(0.000) & 1.000(0.000) &0.998(0.000)  \\
        & &  & 100 & 1.000(0.000) & 0.998(0.000)& 1.000(0.000) \\
        & & \multirow{2}{*}{100} & 50 & 1.000(0.000) & 0.998(0.000) &0.998(0.000)  \\
        & &  & 100 & 1.000(0.000) & 1.000(0.000)& 1.000(0.000) \\
      
        & \multirow{4}{*}{8} & \multirow{2}{*}{50} & 50 & 1.000(0.000) & 1.000(0.000) &1.000(0.000)  \\
        & &  & 100 & 1.000(0.000) & 1.000(0.000)& 1.000(0.000)\\
        & & \multirow{2}{*}{100} & 50 & 1.000(0.000) & 1.000(0.000) &1.000(0.000) \\
        & &  & 100 & 1.000(0.000) & 1.000(0.000)& 1.000(0.000)\\

         \toprule
    \end{tabular}}
    \caption{The frequencies of exact estimation and underestimation of the number of weak factors under Scenario ($\bm{a}$) over 500 replications. }
    \label{tab:weak number}
\end{table}

\section{Empirical Applications}\label{sec:real data}
\subsection{Macroeconomic dataset}\label{subsec:real1}
In this section, we illustrate the empirical usefulness of the TransPCA method by analyzing a U.S. monthly macroeconomic dataset  available from Federal Reserve Economic Data \url{https://fred.stlouisfed.org/}. The dataset contains 126 macroeconomic variables ($N=126$) from January 1959 to August 2024, covering a total of 788 months.  One may refer to the supplement for detailed information about the macroeconomic variables. 

\begin{figure}[!h]
    \centering
    \begin{minipage}[t]{0.7\textwidth}  
        \centering
        \includegraphics[width=1\linewidth]{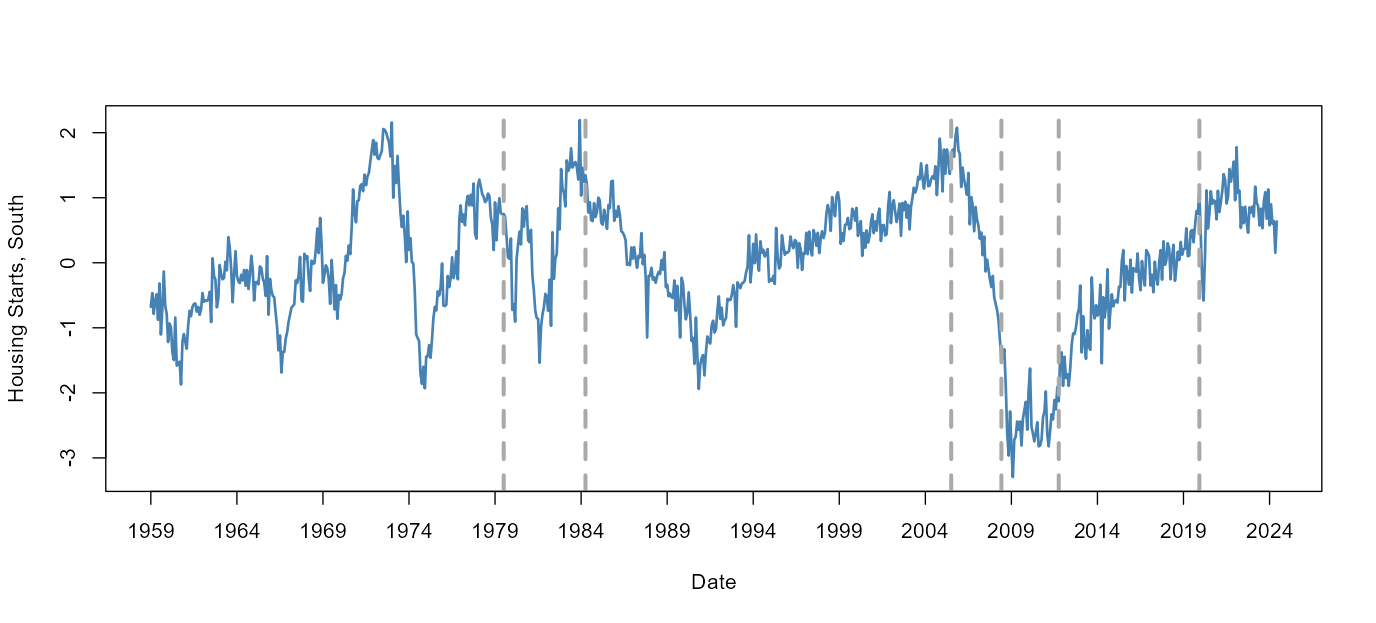} \\
    \end{minipage} \par \vspace{-2em} 
    \begin{minipage}[t]{0.7\textwidth}  
        \centering
        \includegraphics[width=1\linewidth]{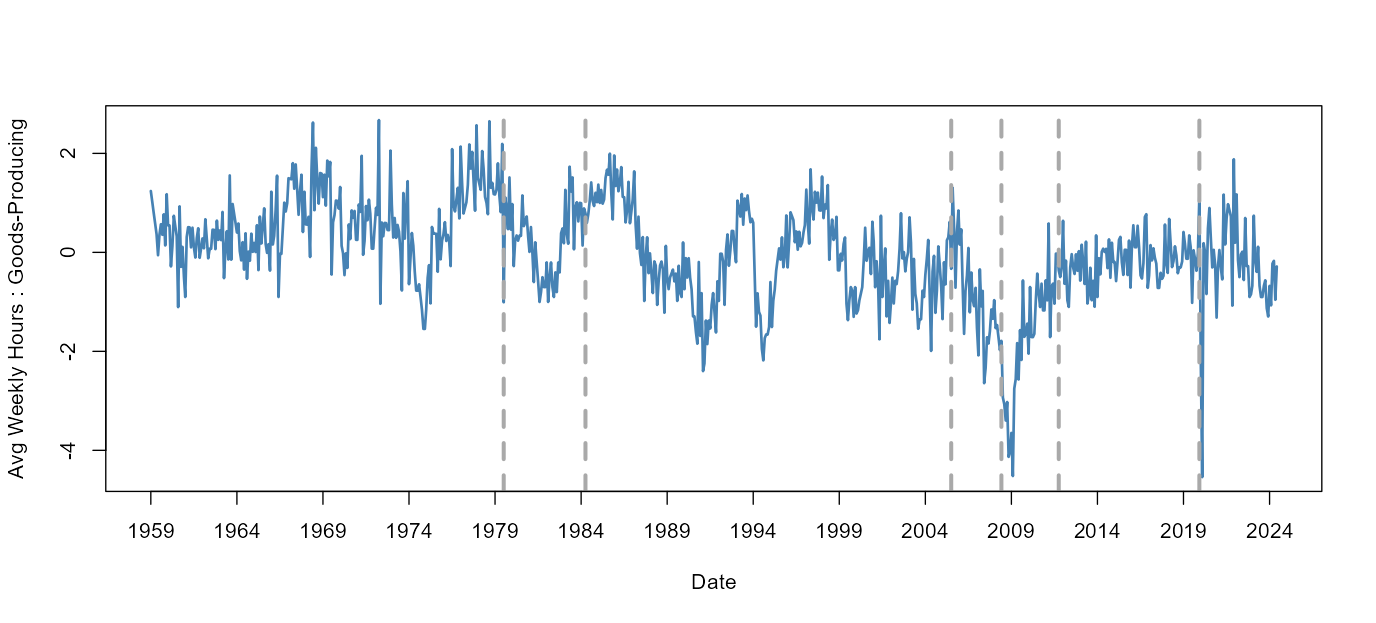} \\
    \end{minipage} \par \vspace{-1em}
    \caption{Time series plots for Housing Starts (South) and Avg Weekly Hours : Goods-Producing. The gray dashed vertical lines are the change point locations. }  
    \label{figure:change points}
    \vspace{-0.5em}
\end{figure}

We first  preprocess the raw dataset, which involves handling non-stationarity, imputing missing values, and dealing with outliers. Given the exceptionally extended time horizon, it is inevitable that the macroeconomic structure would undergo significant structural changes. Therefore we first identify the change points by the method introduced in \cite{Baltagi2017identification}, which gives us 6 change points \textendash \ September 1979, June 1984, September 2005, August 2008, December 2011, February 2020 \textendash \ and  divides the time series into seven segments. The identified change time points includes some well-known significant global events, such as the oil crisis in 1979, the global financial crisis in 2008, and the COVID-19 pandemic in 2020.  Figure \ref{figure:change points} shows the time series of two macroeconomic variables, from which it can be seen  that fluctuations occur at the identified change time points emphasized by gray dashed lines.
The final segment, from March 2020 to August 2024, is viewed  as the target dataset $\Xb^{(0)}\in \RR^{54 \times 126}$ ($T_0=54$), while the remaining segments, ordered chronologically, are treated as the auxiliary datasets $\Xb^{(k)}$, for $k\in [6]$. The Eigenvalue Ratio method in \cite{ahn2013eigenvalue} suggests that $(r_1,r_2,r_3,r_4,r_5,r_6)=(1,1,1,2,3,4)$. For the target dataset, we estimate the factor strength by the estimator  $\hat{\alpha}_i^\ast$ introduced in Section \ref{sec:Estimate FS}. It turns out that for the target dataset $\Xb^{(0)}$, there exist  two weak factors, with factor strengths 0.76 and 0.56 respectively. 

We estimate the target factor model using target-only PCA, blind TransPCA and non-oracle TransPCA as in the simulations, respectively. For the non-oracle TransPCA method, we use 10-fold cross-validation to select parameter $\tau$,  and identify the informative datasets according to Algorithm \ref{alg2}. The results indicate that the informative datasets are $\Xb^{(4)}$, $\Xb^{(5)}$ and $\Xb^{(6)}$. To compare these methods, we employ a rolling-validation procedure. For each month $t$, we repeatedly use the past $n$ (bandwidth) months of observations to estimate the factor loading 
matrix $\hat{\bLambda}_t^n$. The estimated loading matrix is then used to obtain the factors $\hat{\bm{f}}_t$ and corresponding common components $\hat{\bLambda}_t^n\hat{\bm{f}}_t$. We further define the average MSE as  
$$
\text{MSE}_n=\frac{1}{N(T_0-n)}\sum_{t=n+1}^{T_0}\big\|\bX_t-\hat{\bLambda}_t^n\hat{\bm{f}}_t\big\|^2.
$$
Table \ref{tab:9} presents the MSE values of different methods as the bandwidth varies. It can be seen that the MSE of the non-oracle TransPCA method is smaller than those of both the blind TransPCA and target-only methods, with the blind TransPCA method having a larger MSE than the target-only method. This further validates the fact that when  certain level of similarity between auxiliary datasets and target dataset exists,  the TransPCA method yields better performance than performing PCA solely on the target dataset.   Furthermore, when useless datasets are included in the source panels, the non-oracle TransPCA method with dataset selection capability would avoid negative transfer and thus be more effective.
\renewcommand\arraystretch{1.2}
\begin{table}[!h]
    \centering
    \scalebox{1}{
    \begin{tabular}{ccccccccccc}
    \toprule
        $n$ & non-oracle TransPCA&target-only & blind TransPCA  \\
        \hline
        10& \textbf{0.464} & 0.493 & 0.526\\
       15& \textbf{0.446} & 0.480 & 0.514\\
        20& \textbf{0.461} & 0.496 & 0.522\\
         \toprule
    \end{tabular}}
    \caption{The average MSE of different methods for bandwidths $n=10$, 15 and 20.}
    \label{tab:9}
\end{table}



\subsection{Real portfolio returns dataset}
In this section, we analyze a financial  dataset consisting of monthly returns of 100  portfolios. The target dataset consists of monthly average  weighted returns of 100 portfolios, which 
are the intersections of 10 portfolios formed
based on size (market equity, ME) and 10 portfolios formed based on the ratio of book equity to market equity (BE/ME). 
The dataset covers a total of 53 months, from August 2020 to December 2024. Two auxiliary datasets $\Xb^{(1)}$ and $\Xb^{(2)}$ consist of two types of portfolios respectively: one formed by the intersections of 10 size  portfolios and 10 profitability portfolios, covering 180 months from January 2010 to December 2024;
 the other comprised of portfolios formed  by the intersections of 10 size portfolios and 10 investment portfolios, covering 70 months from January 2010 to October 2015. Therefore, for this real financial expample, we have $(T_0,T_1,T_2)=(53,180,70)$ and  $N=100$.


     
        

We first adjust the return series by subtracting the corresponding monthly market
excess returns and then standardizing each of the series. Next, we impute the missing values using the factor-model-based method proposed by \cite{xiong2023large}. The augmented Dickey-Fuller test rejects the null hypotheses for all
the three panels of series, indicating their stationarity. 
{The Eigenvalue-Ratio (ER) method in \cite{ahn2013eigenvalue} indicates that $(r_0,r_1,r_2)=(1,1,1)$ for the three preprocessed datasets, whereas our proposed factor strength estimation method suggests that the target factor strengths $(\hat{\alpha}_1^\ast,\hat{\alpha}_2^\ast, \hat{\alpha}_3^\ast)=(0.81,0.57,0.35)$. For better illustration, we take both $r_0=2$ and $r_0=3$. Regardless of   $r_0=2$ or $r_0=3$, the number of shared ``weak factors'' estimated by our TransED method is $\hat{s}=1$. In other words,  the auxiliary datasets only contribute to the estimation of one weak factor in the target model.} The informative datasets selected by the
non-oracle TransPCA procedure
are $\Xb^{(1)}$ and $\Xb^{(2)}$, meaning that both sources are informative.

\begin{figure*} 
  \centering
  \begin{subfigure}{0.4\linewidth} 
\includegraphics[width=0.98\linewidth]{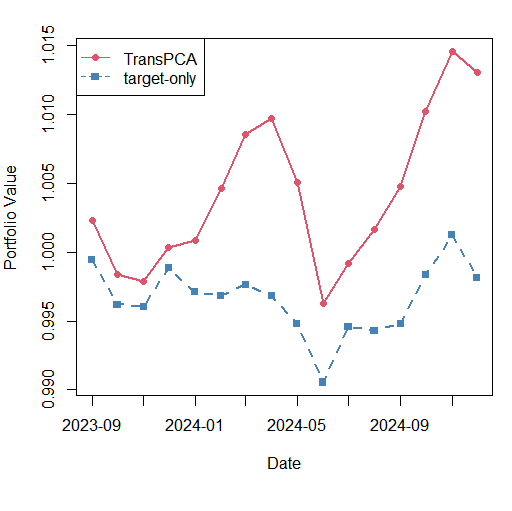} 
    \caption{$r_0=2$}
  \end{subfigure}
  \begin{subfigure}{0.4\linewidth}
\includegraphics[width=0.98\linewidth]{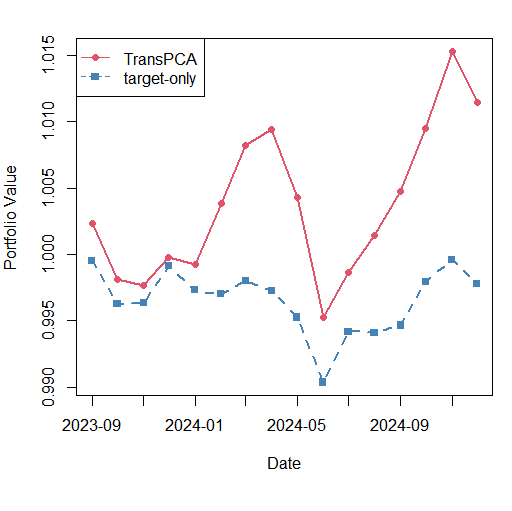} 
    \caption{$r_0=3$}
  \end{subfigure}
  \caption{Portfolio value curves using TransPCA and target-only methods with $r_0=2$ (left) and $r_0=3$ (right).}
  \label{fig:realdata2}
\end{figure*}


We design a rolling investing scheme to compare  the performances of  TransPCA and target-only methods in terms of monthly profitability 
by constructing risk-minimization portfolios.
Given the covariance matrix $\bSigma$ of the  portfolios, a risk-minimization portfolio weights of portfolios is given by
$$
\hat{\bw} = \argmin_{\bm{1}^{\top}\bw=1}\bw^{\top}\bSigma \bw=\frac{\bSigma^{-1}\bm{1}}{\bm{1}^{\top}\bSigma^{-1}\bm{1}}.
$$
However, the covariance matrix $\bSigma$ is unknown in practice and needs to be estimated. To this end, we assume a latent factor structure for the portfolios panel and utilize both the TransPCA and target-only methods to estimate the factor structure, respectively. In detail, for each month $t$, we recursively train the factor model using returns from the past 30 months with both methods. The estimated common components  and idiosyncratic errors are both of dimension $30 \times 100$, denoted as $\hat{\Cb}_t^{(0)}$ and $\hat{\Eb}_t^{(0)}$, respectively.  Then, we estimate the covariance matrix $\hat{\bSigma}_t$ at month $t$ by 
$$
\hat{\bSigma}_t=\frac{1}{30}\hat{\Cb}_t^{(0)\top}\hat{\Cb}_t^{(0)}+\text{HardThresh}\left(\frac{1}{30}\hat{\Eb}_t^{(0)\top}\hat{\Eb}_t^{(0)}\right),
$$
where $\text{HardThresh}\left(\hat{\Eb}_t^{(0)\top}\hat{\Eb}_t^{(0)}/{30}\right)$ represents the hard-threshold operator proposed by \cite{Peter2008Covariance} in view of the sparsity of the covariance matrix for the idiosyncratic errors. The portfolio weights $\hat{\bw}_t$  are calculated with $\hat{\bSigma}_t$ and  the portfolio return for month $t$ is given by $\hat{\bw}_t^{\top}\bX_t$, where $\bX_t$ is the return vector of the portfolios. Figure \ref{fig:realdata2} shows portfolio value curves of the portfolios   by this strategy from September 2023 to December 2024. As can be seen, our TransPCA method, aided by two auxiliary datasets, results in higher portfolio returns than the target-only method, regardless of  $r_0=2$ or $r_0=3$.
This real financial example also supports the necessity  of taking the availability of large
number of auxiliary source datasets into account in factor modeling. These auxiliary finance panels, after screening, would greatly help in target factor analysis. This work then serves as a new complement to all financial econometrics with factor modeling.

\section{Conclusions}\label{sec:conclusion}
In this paper 
we propose a TransPCA method to estimate the target factor model with weak factors
by leveraging useful information from a large body of panel  datasets. We propose a novel model setup in which  the target loading space corresponding to the weak factors is shared or similar to  the loading spaces of the auxiliary datasets,
and integrate useful information across auxiliary datasets via manipulating with a novel weighted average of the  projection matrices, which is the core idea and contribution of this work. 
Theoretical results show that when weak factors exist, the TransPCA estimators outperform the traditional PCA estimators derived solely based on the target dataset under certain assumptions. TransPCA estimators  would also achieve  superior or at least comparable convergence rates with those of PCA estimators in \cite{bai2003inferential} when all factors of the target model are strong. We innovatively propose a TransED method based on the maximum eigenvalue gap of the weighted average projection matrix  to determine the number of weak factors and introduce two different methods for estimating the factor strengths of the target model.  We propose a criterion to select informative auxiliary panels to avoid negative transfer and provide the corresponding theoretical guarantee. 

In the last few years, matrix or tensor time series are growing popular and matrix/tensor factor models have been proposed for achieving multi-mode dimension reduction simultaneously \citep{wang2019factor,chen2023statistical,he2023one,he2024online}. 
In the future, we tend to investigate how to leverage useful information from auxiliary datasets to improve estimation of matrix/tensor factor model of the target dataset. This would be more challenging due to the interplay of multi-modes and warrants further research.

\section*{Acknowledgements}
This work is supported by NSF China (12171282) and Qilu Young Scholars Program of Shandong University, China.

\bibliographystyle{model2-names}
\bibliography{ref}

\end{document}